\documentclass[a4paper,12pt]{amsart}
\usepackage{amscd,amsthm,amsfonts,latexsym,amssymb,amsmath}

\newtheorem{thm}{Theorem}[section]
\newtheorem{lem}[thm]{Lemma}
\newtheorem{cor}[thm]{Corollary}
\newtheorem{prop}[thm]{Proposition}

\theoremstyle{definition}

\newtheorem{defn}[thm]{Definition}
\newtheorem{rem}[thm]{Remark}
\newtheorem{exm}[thm]{Example}
\newtheorem{notn}[thm]{Notation}

\numberwithin{equation}{section}

\newcommand{\vvec}[1]{\underline{\boldsymbol{#1}}}	

\newcommand{\dd}{{\mathrm d}}
\newcommand{\pa}{{\partial}}

\newcommand{\zbar}{{\bar{z}}} 

\newcommand{\C}{{\mathbb C}}  

\newcommand{\CP}{{\mathbb C}P} 

\newcommand{\U}{\mbox{\rm U}}   

\newcommand{\uu}{{\mathfrak u}}
	
\newcommand{\gl}{\mathbf{gl}}  

 \newcommand{\half}{\textstyle\frac{1}{2}}  

\newcommand{\ul}{\underline}

\newcommand{\CC}{\underline{\mathbb C}} 

\newcommand{\End}{\text{\rm End}}  
\newcommand{\rk}{\text{\rm rank}\,}  

\newcommand{\HH}{\mathcal H}

\renewcommand{\phi}{\varphi}

\newcommand{\la}{\lambda}

\newcommand{\al}{\alpha}

\newcommand{\ga}{\gamma}
\newcommand{\Ga}{\Gamma}

\newcommand{\si}{\sigma}

\newcommand{\Om}{\Omega}

\newcommand{\aal}{\underline{\alpha}} 
\newcommand{\pphi}{\underline{\phi}} 
\newcommand{\hh}{\underline{h}} 
\newcommand{\BB}{\underline{B}} 

\newcommand{\p}{{\perp}} 

\newcommand{\wt}{\widetilde} 
\newcommand{\ov}{\overline} 

\DeclareMathOperator{\spn}{\mbox{\rm span}} 
\DeclareMathOperator{\iim}{\mbox{\underline{\rm Im}}} 

\DeclareMathOperator{\kker}{\mbox{\underline{\rm ker}}}

\newcommand{\XX}{\underline{X}} 
\newcommand{\WW}{\underline{W}} 

\newcommand{\M}{\mathcal M} 

\newcommand{\alg}{\mbox{\tiny\rm alg}} 

\begin{document}

\title[All harmonic $2$-spheres in the unitary group]{All harmonic $2$-spheres in the unitary group, completely explicitly}

\author{Maria Jo\~ao Ferreira}
\author{Bruno Ascenso Sim\~oes}
\author{John C.\ Wood}

\address{\textit{MJF}, \textit{BAS}: Centro de
Matem\'atica e Aplica\c{c}\~oes Fundamentais,
Universidade de Lisboa, Av.\ Prof.\ Gama Pinto 2\\1649-003, Lisbon, Portugal}

\email{mjferr@ptmat.fc.ul.pt; b.simoes@ptmat.fc.ul.pt}

\address{\textit{JCW}: Department of Pure Mathematics, University of Leeds\\
Leeds LS2 9JT, Great Britain}

\email{j.c.wood@leeds.ac.uk}

\thanks{BAS thanks the \textit{Funda\c{c}\~ao para a Ci\^encia e Tecnologia}, Portugal, for financial support.
JCW thanks the Gulbenkian foundation, \textit{CMAF}
 and the Faculdade de Ci\^encias, Universidade de Lisboa,
and \textit{IMADA}, University of Southern Denmark, Odense, for support and hospitality during the preparation of this work.}

\maketitle

\begin{abstract}
We give a completely explicit formula for all harmonic maps of finite uniton number from a Riemann surface to the unitary group $\U(n)$ in any dimension, and so all harmonic maps from the $2$-sphere, in terms of freely chosen meromorphic functions on the surface and their derivatives, using only combinations of projections and avoiding the usual $\ov{\pa}$-problems or loop group factorizations.  We interpret our constructions using Segal's Grassmannian model, giving an explicit factorization of the algebraic loop group, and showing how to obtain harmonic maps into a Grassmannian.

\end{abstract}

\section*{Introduction}
In \cite{Uhl}, K.\ Uhlenbeck showed how to construct all harmonic maps from the $2$-sphere to the unitary group
$\U(n)$, equipped with its standard bi-invariant metric,
by starting with a constant map and modifying it by successively multiplying by suitable maps into Grassmannians --- a process called \emph{adding a uniton} or, by others, \emph{flag transform}.  The unitons must be holomorphic with respect to the Koszul--Malgrange holomorphic structure coming from the current harmonic map; in previous papers, for example, \cite{HeSh,Wo-Un}, these successive holomorphic structures were found in terms of the previous one and the solution to a $\ov{\pa}$-problem to which, in general, no explicit solution could be given.  We solve this problem by imposing the \emph{covering condition} of G.\ Segal \cite{Seg} on each uniton, which says that it must project surjectively onto the previous one; we can then find the next holomorphic structure explicitly just by using suitable projections of meromorphic functions.  We can thus build harmonic maps explicitly in terms of freely chosen meromorphic functions and their derivatives, see
Theorem \ref{thm:constr-H}.  We then give some examples, and a reformulation
where the data becomes a holomorphic subbundle (Theorem \ref{thm:constr-B}) related to the Grassmannian model.  Note that, by thinking of them as stationary Ward solitons, B.\ Dai and C.-L. Terng \cite{DaTe} obtain an explicit formula for the unitons of the Uhlenbeck factorization, but although the two factorizations are dual (see Remark \ref{rem:dual}), the formulae do not seem to be equivalent under that duality.

In Section \ref{sec:unitons}, after recalling the theory of Uhlenbeck, we explain how any possible factorization by covering unitons is given by our data, see Theorem \ref{thm:all-covering}.

To prove that our formulae give \emph{all} harmonic maps from the
$2$-sphere, we must show how to factorize an arbitrary harmonic map into the product of unitons satisfying the covering condition.
Uhlenbeck's factorization uses the kernel of the constant term $T_0$ of the extended solution as the next uniton.  We employ the dual of this factorization which uses the image of the adjoint of the highest coefficient $T_r$ of the extended solution, instead. This factorization was studied by \cite{HeSh} where the unitons were called \emph{AUN-flag factors} and were characterized in terms of $T_r$; our covering condition is equivalent, but depends only on the successive unitons.  As for Uhlenbeck's factorization, the covering condition ensures
\emph{uniqueness} of the factorization, see Theorem \ref{thm:unique-fact},
cf.\ \cite{Seg}.

A harmonic map from an arbitrary Riemann surface is said to be \emph{of finite uniton number} \cite{Uhl} if there exists a factorization into (a finite number of) unitons; our work applies equally well to such maps.

In \cite{Seg}, G.\ Segal gave a model for the loop group of $\U(n)$ as an (infinite-dimensional) Grassmannian, and showed that harmonic maps of finite uniton number correspond to holomorphic maps into a related finite-dimensional Grassmannian.
In Section \ref{sec:loop}, we calculate the Grassmannian model solution corresponding to our harmonic maps in the form discussed by
M.\ Guest \cite{Gu-update}, and give an interpretation
in terms of an explicit Iwasawa factorization of the algebraic loop group.

In \cite{BuGu}, F.E. Burstall and M. Guest gave a construction of harmonic $2$-spheres into any compact Lie group in terms of holomorphic data.
In \cite{Gu-update}, M. Guest interpreted that work in terms of the Grassmannian model; we compare our formulae with his.

Our construction includes all harmonic maps into complex Grassmannians; at the end of Section \ref{sec:loop} we see how to obtain these via the Grassmannian model. A more geometrical approach will be the subject of a future paper.

To make the paper more accessible, the tools needed are introduced only when they are required.  Only elementary concepts are required to understand the statements of the main theorems in Section \ref{sec:build}.  To explain where our formulae come from,  we discuss \emph{unitons} in Section
\ref{sec:unitons}, but without using extended solutions.  In Section
\ref{sec:fact}, we need \emph{extended solutions} to show that all harmonic maps are obtained.  Only in Section \ref{sec:loop} is the \emph{loop group} and its \emph{Grassmannian model} introduced.

In \cite{SvWo} M.~Svensson and the third author
further interpret and develop the work in this paper to obtain explicit
formulae for all harmonic $2$-spheres in the symplectic and
special orthogonal groups, and their inner symmetric spaces.
The authors are grateful to Martin Svensson for some very useful discussions on this work,  to Bo Dai for illuminating correspondence on \cite{DaTe}, and to the referee for some pertinent comments.

\section{Explicit formulae for harmonic maps}\label{sec:build}

Let $M^2$ be a Riemann surface.   Let $n \in \{1,2,\dots\}$, and let $\CC^n$ denote the trivial complex bundle $M^2 \times \C^n$ equipped with the standard Hermitian inner product: $\langle \vec{u},\vec{v} \rangle = u_1 \ov{v}_1 + \cdots + u_n \ov{v}_n$ \ $\bigl(\vec{u} = (u_1,\ldots, u_n), \vec{v} = (v_1,\ldots, v_n) \in \C^n\bigr)$ on each fibre.    For a subbundle $\ul{\al}$ of $\CC^n$, let $\pi_{\al}$ and $\pi_{\al}^{\p}$ denote orthogonal projection onto $\aal$ and onto its orthogonal complement $\aal^{\p}$, respectively. By a \emph{$\C^n$-valued meromorphic function} or \emph{meromorphic vector} $H$ on $M^2$, we simply mean an $n$-tuple of meromorphic functions; we denote its $k$'th derivative with respect to some local complex coordinate on $M^2$ by $H^{(k)}$.
With explanations to follow, we can state our explicit construction of harmonic maps as follows.

\begin{thm} \label{thm:constr-H} For any $r \in \{0,1,\ldots, n-1\}$, let $(H_{i,j})_{0 \leq i \leq r-1,\ 1 \leq j \leq n}$ be an $r \times n$ array of\/ $\C^n$-valued meromorphic functions on $M^2$, and let $\phi_0$ be an element of\/ $\U(n)$.
For each $i = 0,1, \ldots, r-1$, set $\aal_{i+1}$ equal to the subbundle of\/ $\CC^n$ spanned by the vectors
\begin{equation} \label{alpha-H}
K^{(k)}_{i,j} = \sum_{s=k}^{i} C^{i}_s H^{(k)}_{s-k, j} \qquad (j = 1,\ldots, n, \
 k =0,1, \ldots, i).
\end{equation}

Then the map  $\phi:M^2 \to \U(n)$ defined by
\begin{equation} \label{fact}
\phi = \phi_0 (\pi_1 - \pi_1^{\p}) \cdots (\pi_r - \pi_r^{\p})
\end{equation}
is harmonic.

Further, all harmonic maps from $M^2$ to $\U(n)$ of finite uniton number, and so \emph{all} harmonic maps from $S^2$ to $\U(n)$, are obtained this way.
\end{thm}

\begin{notn} \label{notn:notn}
(i) We write $\pi_i = \pi_{\al_i}$ and
$\pi_i^{\p} = \pi_{\al_i}^{\p}$.

(ii) For integers $i$ and $s$ with $0 \leq s \leq i$, $C^i_s$ denotes the \emph{$s$'th elementary function} of the projections $\pi_i^{\p}, \ldots, \pi_1^{\p}$ given by
\begin{equation} \label{Cis}
C^i_s = \sum_{1 \leq i_1 < \cdots < i_s \leq i}
	\pi_{i_s}^{\p} \cdots \pi_{i_1}^{\p}\,.
\end{equation}
For example,
$C^2_1 = \pi_1^{\p} + \pi_2^{\p}$ and
$C^3_2 = \pi_2^{\p}\pi_1^{\p} + \pi_3^{\p}\pi_1^{\p} +\pi_3^{\p}\pi_2^{\p}$.  We define $C^i_s$ to be the identity when $s = 0$ and zero when $s <0$ or $s > i$.  Note that the $C^i_s$ satisfy a property like that for Pascal's triangle:
\begin{equation}\label{C-Pascal}
C^i_s = \pi_i^{\p} C^{i-1}_{s-1}+ C^{i-1}_s  \quad
(i \geq 1, \ 0 \leq s \leq i).
\end{equation}

(iii) For any holomorphic bundle $Z$ over $M^2$, we shall write $\M(Z)$ for the space of its meromorphic sections, a vector space over the field of meromorphic functions on $M^2$.   Note that $Z$ can be recovered as the span of a basis for that vector space.

More generally, given a finite collection $C$ of meromorphic sections of a holomorphic bundle $Z$, away from a discrete set $D$ which includes the poles and zeros of the meromorphic vectors, the span of $C$ has constant rank and so defines a subbundle of $Z$ over $M \setminus D$.  By `filling in holes' as in \cite[Proposition 2.2]{BuWo}, such a subbundle can be extended smoothly over $D$ to a bundle over the whole of $M^2$; we denote this by $\spn C$.  The sections $K^{(k)}_{i,j}$ above are meromorphic sections of the trivial bundle with respect to a suitable holomorphic structure, see
Proposition \ref{prop:calcs} below.

(iv) For $i = 0,1,2,\ldots$ and $k=0,\,\ldots, i$, we let $\aal_{i+1}^{(k)}$ denote the subbundle  $\spn\{K^{(k)}_{i,j} : j=1,\ldots, n\}$; then $\aal_{i+1}$ is the sum
$\aal_{i+1} = \sum_{k=0}^{i}\aal_{i+1}^{(k)}$\,.

We shall see (Remark \ref{rem:gen-uniton}) that the $\aal_{i+1}^{(k)}$
are independent of choice of complex coordinate used to compute
the derivatives.
\end{notn}

\begin{exm} (\emph{The first three subbundles})

The first subbundle $\aal_1$ is formed from the first row of the data, precisely,
\begin{equation} \label{alpha1}
\aal_1 = \aal_1^{(0)} = \spn\{H_{0,j}\}\,.
\end{equation}

The second subbundle $\aal_2$ is formed from the first two rows of the data and their first derivatives; precisely,
  $\aal_2 = \aal_2^{(0)} + \aal_2^{(1)}$ with
\begin{align}
\aal_2^{(0)} &= \spn\{H_{0,j} + \pi_1^{\p}H_{1,j} \}\,, \label{alpha20}\\
\aal_2^{(1)} &= \spn\{\pi_1^{\p}H_{0,j}^{(1)} \}\,. \label{alpha21}
\end{align}

The third subbundle $\aal_3$ is formed from the first three rows of the data and their first and second derivatives; precisely,
  $\aal_3 = \aal_3^{(0)} + \aal_3^{(1)} + \aal_3^{(2)}$ with
\begin{align}
\aal_3^{(0)} &= \spn\{H_{0,j} + (\pi_1^{\p} + \pi_2^{\p})H_{1,j}
	+ \pi_2^{\p}\pi_1^{\p}H_{2,j} \}\,,
\label{alpha30} \\
\aal_3^{(1)} &= \spn\{(\pi_1^{\p}+\pi_2^{\p})H_{0,j}^{(1)}
	+ \pi_2^{\p}\pi_1^{\p}H_{1,j}^{(1)} \}\,, \label{alpha31}\\
\aal_3^{(2)} &= \spn\{\pi_2^{\p}\pi_1^{\p}H_{0,j}^{(2)}\}\,. \label{alpha32}
\end{align}
\end{exm}

\begin{exm} (\emph{Maps into Grassmannians}) \label{exm:Grass}
For any $k \in \{0,1,\ldots, n\}$, let $G_k(\C^n)$
denote the \emph{complex Grassmannian} of $k$-dimensional subspaces of $\C^n$ equipped with its standard structure as a Hermitian symmetric space.  Then there is a one-to-one correspondence between smooth subbundles $\aal$ of $\CC^n$ of rank $k$ and smooth maps $\al:M^2 \to G_k(\C^n)$ given by setting the fibre of $\aal$ at $p \in M^2$ equal to the
subspace $\al(p)$; the subbundle is holomorphic
if and only if the map is.   It is convenient to denote the disjoint union $\cup_{k=0}^n G_k(\C^n)$ by $G_*(\C^n)$.

The \emph{Cartan embedding}
$\iota:G_*(\C^n) \to \U(n)$ is given by
$\iota(A) = \pi_A - \pi_A^{\p}$.
It is isometric and totally geodesic, so that \emph{a smooth map $\phi$ into a Grassmannian is harmonic if and only
$\iota \circ \phi = \pi_{\phi} - \pi_{\phi}^{\p}$
is harmonic into $\U(n)$}; note that $\phi$ is harmonic if and only if $\phi^{\p}$ is harmonic and
$\iota \circ \phi^{\p} = -\iota \circ \phi$.  In particular, our formulae give all harmonic maps into Grassmannians, see \S \ref{sec:loop}B.  We shall often consider a map into a Grassmannian $G_*(\C^n)$ as a map into $\U(n)$ via the Cartan embedding, without comment.
\end{exm}

Say that two harmonic maps $\phi,\wt{\phi}:M^2 \to \U(n)$ are
\emph{left-equivalent} if $\phi = C \wt{\phi}$ for some constant
$C \in \U(n)$.  We shall classify harmonic maps up to left-equivalence, thus we can assume that $\phi_0 = I$ in \eqref{fact}.

\begin{exm} \label{exm:constr-low}
(i) \emph{The case $r=0$.}  There is no data $H_{i,j}$ and $\phi$ is a constant map.

\smallskip

(ii) \emph{The case $r=1$.}  The array $(H_{i,j})$ consists of a single row $(H_{0,j})$ of meromorphic vectors; $\aal_1$ is the holomorphic subbundle of $\CC^n$ spanned by these vectors.
Writing $d_1 = \rk\aal_1$, the harmonic map $\phi$ is equal to the holomorphic map $\al_1:M^2 \to G_{d_1}(\C^n)$.

Any non-constant harmonic map $S^2 \to \U(2)$ is of this form
with $d_1=1$ (see Example \ref{exm:low-dim}).

\smallskip

(iii) \emph{The case $r=2$.}
The data $(H_{i,j})$ consists of $2$ rows; $\aal_1$ and $\aal_2$ are given by
\eqref{alpha1}--\eqref{alpha21}.
Note that $\aal_2^{(1)} = \spn\{\pi_1^{\p}H_{0,j}^{(1)}\}$
 is the \emph{($\pa'$-)Gauss bundle}  $G^{(1)}(\aal_1)$
of $\aal_1$ defined by
$G^{(1)}(\aal_1) = \{\pi_1^{\p} H^{(1)} : H \in \M(\aal_1)\}$.

In the special case that all $H_{1,j}$ are zero, then
$\aal_1 \subseteq \aal_2$; in fact,
$$
\aal_2 = \aal_1 \oplus G^{(1)}(\aal_1)
= \spn\{\si, \si^{(1)}: \si \in \M(\aal_1)\}
$$
is the \emph{first associated curve of\/ $\aal_1$},  and
$\pphi_2 = \aal_1^{\p} \cap \aal_2 = G^{(1)}(\al_1)$.

The cases $r \leq 2$ account for all harmonic maps from $S^2$ to $\U(3)$. By allowing also $r=3$ we obtain all harmonic maps $S^2 \to \U(4)$, see Example \ref{exm:low-dim} for precise formulae.

\smallskip

(iv) \emph{General $r$, one non-zero row.}  If the only non-zero entries
in the array $(H_{i,j})$ are  $H_{0,1},\ldots, H_{0,d_1}$, then $\aal_1$ is the holomorphic subbundle spanned by these vectors and
each $\aal_{i+1}$ is the \emph{$i$'th associated curve of\/ $\aal_1$}.   Here, the $i$'th associated curve of a holomorphic subbundle $\hh$ of $\CC^n$ is defined by
\begin{equation} \label{assoc-i}
\hh_{(i)} = \spn\{\si^{(j)}: j=0,1,\ldots, i,\ \si \in \M(\hh)\}\,.
\end{equation}
The resulting harmonic map $\phi$ is a harmonic map into a Grassmannian described as follows. For any holomorphic subbundle $\hh$ of
$\CC^n$, let $G^{(i)}(\hh)$ denote the
\emph{$i$'th Gauss bundle of $\hh$} given by
\begin{equation} \label{Gauss-i}
G^{(i)}(\hh) = {\hh_{(i-1)}}^{\!\!\p} \cap \hh_{(i)} \,.
\end{equation}
Then the formula
$\sum_{i=0}^{[(r-1)/2]}G^{(r-2i-1)}(\aal_1)$
gives $\pphi$ if $r$ is odd or $\pphi^{\p}$ if $r$ is even.
See also Example \ref{exm:isotropic}.
\end{exm}

\begin{rem} \label{rem:constr-H}
(i) The subbundles $\aal_i$\,, or the corresponding
 factors $\pi_i - \pi_i^{\p}$ in \eqref{fact},
are \emph{unitons} and \eqref{fact} is a \emph{uniton factorization},
see Section \ref{sec:unitons}.

(ii) We are not insisting that any of the functions $H_{i,j}$ should be non-zero, so it could happen that some $\aal_i$ are the zero subbundle.  It could also happen that some $\aal_i$ are the whole of $\CC^n$; in both cases, the factor $\pi_i - \pi_i^{\p}$ in \eqref{fact} could be removed.
On excluding this case and imposing a covering and fullness condition, the unitons $\aal_i$ \emph{are} uniquely determined by the map $\phi$, see Theorem \ref{thm:unique-fact}.

(iii) We call the subbundle $\aal_{i+1}^{(0)}$ the \emph{generating subbundle
of\/} $\aal_{i+1}$; it is spanned by\\[-4ex]
\begin{equation} \label{K}
K_{i,j} := K_{i,j}^{(0)} = \sum_{s=0}^{i} C^{i}_s H_{s,j}\,.
\end{equation}
The other subbundles $\aal_{i+1}^{(k)}$ are determined by it; also, the possible subbundles $\aal_{i+1}^{(0)}$ for a given $\aal_{i+1}$ can be characterized, see Remark \ref{rem:gen-uniton}.

(iv)   The array $(H_{i,j})$ which determines a given list of
$\aal_i$ is not unique; for example, it can be replaced by any array with the same column span over the meromorphic functions.  Indeed, by column operations one can replace the array by one with linearly independent columns in the `echelon' form:
\begin{equation*}
\left[
\begin{array}{ccccccccccc}
H_{0,1} & \cdots & H_{0,d_1} & 0 &\cdots & 0
	&0 & \cdots & 0 & 0 & \cdots \\
H_{1,1} & \cdots & H_{1,d_1} & H_{1,d_1+1} &
	\cdots & H_{1,d_2} & 0 & \cdots & 0 & 0 & \cdots \\
H_{2,1} & \cdots & H_{2,d_1} & H_{2,d_1+1} &
	\cdots & H_{2,d_2} & H_{2,d_2+1} & \cdots & H_{2,d_3} & 0 & \cdots \\
\vdots & \vdots & \vdots & \vdots & \vdots & \vdots
	& \vdots & \vdots & \vdots & \vdots & \ddots \\	
\end{array} \right]
\end{equation*}
with $0 \leq d_1 \leq d_2 \leq \cdots \leq d_r \leq n$,
where we further insist that, for each each $i = 1,\ldots, r$, the sub-array made up of the first $i$ rows and first $d_i$ columns has linearly independent columns, equivalently, for each $i$, the vectors $H_{i,d_{i}+1} \ldots, H_{i,d_{i+1}}$ are linearly independent.
Note that $\aal_i^{(0)}$ is constructed from that sub-array by formula \eqref{alpha-H} and so has rank at most $d_i$.

More serious changes can be made; for example, for any given column $(H_0,H_1,\ldots, H_{r-1})^T$, by Lemma \ref{lem:sum=Hi}, we may include the column $(0,H_0,\ldots, H_{r-2})^T$ without changing any $\aal_i^{(k)}$, and by Proposition \ref{prop:calcs}, we may include the column $(0,H_0^{(1)},\ldots, H^{(1)}_{r-2})^T$ without changing any $\aal_i$.  By repeatedly adjoining all columns of the second type, we can make $\aal_i = \aal_i^{(0)}$; see also Remark \ref{rem:gen-uniton}.
\end{rem}

For a more invariant formulation, give $\CC^n$ its standard holomorphic structure so that its meromorphic sections are precisely the meromorphic $\CC^n$-valued functions, and consider $\CC^{rn}$ to be the direct sum of $r$ copies of the holomorphic bundle $\CC^n$. Correspondingly write any $H \in \CC^{rn}$ as  $H= (H_0, H_1, \ldots, H_{r-1})$ with each $H_i \in \CC^n$.
Then the following is equivalent to our main result, Theorem \ref{thm:constr-H}.

\begin{thm} \label{thm:constr-B} For any $r \in \{0,1,\ldots, n-1\}$, let
$\BB$ be a holomorphic subbundle of\/ $\CC^{rn}$, and let $\phi_0$ be an element of\/ $\U(n)$.
For $i = 0,1,\ldots, r-1$, $k = 0,1\ldots, i$, define subbundles\/ $\aal_{i+1}^{(k)}$, $\aal_{i+1}$ by
\begin{equation} \label{alpha-B}
\aal_{i+1}^{(k)} = \Bigl\{\sum_{s=k}^{i} C^{i}_s H^{(k)}_{s-k} : H \in \M(\BB) \Bigr\}, \quad
	\aal_{i+1} = \sum_{k=0}^{i}\aal_{i+1}^{(k)}\,.
\end{equation}
Then the map  $\phi:M^2 \to \U(n)$ defined by \eqref{fact} is harmonic.

Further, all harmonic maps from $M^2$ to $\U(n)$ of finite uniton number, and so \emph{all} harmonic maps from $S^2$ to $\U(n)$, are obtained this way.
\end{thm}

\begin{rem} \label{rem:B-to-H}
The relationship between the two formulations in Theorems \ref{thm:constr-H} and \ref{thm:constr-B} is that, given an array $(H_{i,j})$, $\BB$ is the subbundle spanned by its columns; conversely, given $\BB$, we can choose $(H_{i,j})$ to be any array whose columns span $\BB$.  The subbundle $\BB$ gives the Grassmannian model of $\phi$, see Theorem \ref{thm:Grass-model}.
\end{rem}

\section{Building harmonic maps from unitons}
\label{sec:unitons}

A. \emph{Factorizations by covering unitons}

\medskip

We now explain the origin of our formulae \eqref{alpha-H}, \eqref{alpha-B}.  For this we recall the basic theory of Uhlenbeck, see \cite{Uhl,Wo-Un,Gu-bk}, of harmonic maps from a Riemann surface $M^2$ to the unitary group $\U(n)$.

Let $\phi:M^2 \to \U(n)$ be a smooth map.
Let $\uu(n)$ denote the Lie algebra of $\U(n)$ consisting of the
$n \times n$ skew-Hermitian matrices.
 Define a $1$-form with values in $\uu(n)$ by $A = A^{\phi} = \half \phi^{-1}\dd\phi$\,; thus $A^{\phi}$ is one-half of the pull-back of the (left) Maurer--Cartan form on $\U(n)$.  For convenience, we choose a local complex coordinate $z$ on an open subset of $M^2$; our theory will be independent of that choice. Then we can decompose $A^{\phi}$ into $(1,0)$- and $(0,1)$-parts:
$A^{\phi} = A^{\phi}_z \dd z+ A^{\phi}_{\zbar}\dd\zbar$; note that $A^{\phi}_z$ and  $A^{\phi}_{\zbar}$ are local sections of the endomorphism bundle $\End(\CC^n)$, and each is minus the adjoint of the other.  Set $D^{\phi} = \dd+A^{\phi}$.  Then $D^{\phi}$ is a unitary connection on the trivial bundle $\CC^n$; in fact, it is the pull-back of the Levi-Civita connection on $\U(n)$.

  We write $D^{\phi}_z =\pa_z+A^{\phi}_z$ and $D^{\phi}_{\zbar} =  \pa_{\zbar} + A^{\phi}_{\zbar}$
where $\pa_z = \pa/\pa z$ and $\pa_{\zbar} = \pa/\pa\zbar$.  Give $\CC^n$ the Koszul--Malgrange complex structure \cite{KoMa}; this is the unique holomorphic structure such that a (local) section $\si$ of $\CC^n$ is holomorphic if and only if $D^{\phi}_{\zbar}\si = 0$ for any complex coordinate $z$; we shall denote the resulting holomorphic bundle by $(\CC^n, D^{\phi}_{\zbar})$.  Note that, when $\phi$ is constant, $A^{\phi} = 0$, and the Koszul-Malgrange holomorphic structure is the standard holomorphic structure on $\CC^n$, i.e.,
$(\CC^n, D^{\phi}_{\zbar}) =  (\CC^n, \pa_{\zbar})$.

Since $A^{\phi}_z$ represents the derivative $\pa\phi/\pa z$,
\emph{the map $\phi$ is harmonic if and only if the
 endomorphism $A^{\phi}_z$ is holomorphic}, i.e.,
\begin{equation} \label{hol}
A^{\phi}_z \circ D^{\phi}_{\zbar}= D^{\phi}_{\zbar} \circ A^{\phi}_z\,.
\end{equation}

Let $\phi:M^2 \to \U(n)$ be harmonic and let $\aal$ be a smooth subbundle of $\CC^n$.   Say that $\aal$ is a \emph{uniton} or \emph{flag factor
for $\phi$\/} if
\begin{equation} \label{uniton}
\left\{
\begin{array}{rrl}
{\rm (i)} & D^{\phi}_{\zbar}(\si) \in \Ga(\aal)
	&\mbox{ for all } \si \in \Ga(\aal)\,, \\[0.5ex]
{\rm (ii)} & A^{\phi}_z(\si) \in \Ga(\aal)
	&\mbox{ for all } \si \in \Ga(\aal)\,;
\end{array} \right. \end{equation}
here $\Ga(\cdot)$ denotes the space of smooth sections of a bundle.  These equations say that $\aal$ is a holomorphic subbundle of $(\CC^n, D^{\phi}_{\zbar})$ which is closed under the endomorphism $A^{\phi}_z$.
Uhlenbeck shows \cite{Uhl} that \emph{if\/ $\phi$ is harmonic and $\aal$ is a uniton for $\phi$, then the map $\wt{\phi}:M^2 \to \U(n)$ given by
$\wt{\phi} = \phi(\pi_{\al} - \pi_{\al}^{\p})$ is harmonic}.

  Note that $\aal$ is a uniton for $\phi$ if and only if $\aal^{\p}$ is a uniton for $\wt{\phi}$\,; further
$\phi = - \wt{\phi}(\pi_{\al}^{\p} - \pi_{\al})$\,, i.e., the flag transforms defined by $\aal$ and $\aal^{\p}$ are inverse up to sign.  Furthermore, the connections induced by $\phi$ and $\wt{\phi}$ are related by the simple formulae \cite{Uhl}:
\begin{equation}\label{uniton-DA}
\text{(i)} \quad
A_z^{\wt{\phi}} = A_z^{\phi} + \pa_z\pi_{\al}^{\p}\,,
\quad \text{(ii)} \quad
D_{\zbar}^{\wt{\phi}} = D_{\zbar}^{\phi} - \pa_{\zbar}\pi_{\al}^{\p} \,,	
\end{equation}
which lead to the useful equations:
\begin{equation} \label{commute}
\text{(i)} \quad A_z^{\wt{\phi}}\pi_{\al}^{\p}
	= \pi_{\al}^{\p} A_z^{\phi}\,, \quad
\text{(ii)} \quad \ D_{\zbar}^{\wt{\phi}}\pi_{\al}^{\p}
	= \pi_{\al}^{\p} D_{\zbar}^{\phi} \,.
\end{equation}

As explained in Notation \ref{notn:notn}(iii) and in \cite{Uhl}, the ranks of the kernel and image of the endomorphism $A^{\phi}_z$ on a fibre $\{p\} \times \C^n \subset \CC^n$ are constant as $p$ varies over $M^2$, except on a discrete set where the rank may jump; however, because of the holomorphicity of $A^{\phi}_z$\,, we may extend them over these points to subbundles $\kker A^{\phi}_z$ and $\iim A^{\phi}_z$.

We thus obtain two fundamental examples of unitons \cite{Uhl}:

(i) the \emph{kernel bundle} $\kker A^{\phi}_z$, or any holomorphic subbundle of $(\CC^n,D^{\phi}_{\zbar})$ contained in that bundle; following \cite{Wo-Un}, we call such unitons \emph{basic};

(ii) the \emph{image bundle} $\iim A^{\phi}_z$, or any holomorphic subbundle of $(\CC^n,D^{\phi}_{\zbar})$ containing that bundle; following \cite{PiZa}, we call such unitons \emph{antibasic}.

{}From  \eqref{commute} (cf.\ \cite{Wo-Un}) we see that  \emph{$\aal$ is a basic (resp.\ antibasic) uniton for $\phi$ if and only if\/ $\aal^{\p}$ is an antibasic (resp.\ basic) uniton for $\wt{\phi}$}.

The idea of Uhlenbeck was to start with a constant map $\phi_0:M^2 \to \U(n)$ and build more complicated harmonic maps $\phi_i$ by successively setting
$\phi_i = \phi_{i-1} (\pi_i - \pi_i^{\p})$ \
$(i =1,2,\ldots)$.  Here, as usual, we write $\pi_i = \pi_{\al_i}$ and
$\pi_i^{\p} = \pi_{\al_i}^{\p}$ where $\aal_1,\ldots,\aal_i$ is \emph{a sequence of unitons}, by which we mean that $\aal_{\ell}$ is a uniton
for $\phi_{\ell-1}$ \ $(\ell=1,\ldots, i)$\,.  Thus each $\phi_i$ is a
product of unitons:
\begin{equation} \label{phi_i}
\phi_i = \phi_0 (\pi_1 - \pi_1^{\p}) \cdots (\pi_i - \pi_i^{\p})\,.
\end{equation}

A harmonic map obtained in this way is said to be \emph{of finite uniton number}; Uhlenbeck showed that all harmonic maps from the $2$-sphere to $\U(n)$ are of finite uniton number.

We study sequences of unitons $\aal_1,\ldots,\aal_i$
satisfying the following condition which we call the \emph{covering condition}:
\begin{equation}\label{cover}
\pi_{\ell-1}\aal_{\ell} = \aal_{\ell-1}  \qquad (\ell = 2,\ldots, i)\,.
\end{equation}

 \begin{prop} \label{prop:cover-conds} Let $\phi_i$ be given by
\eqref{phi_i} for some constant map $\phi_0$ and sequence
$\aal_1, \ldots, \aal_i$ of unitons.

Then \eqref{cover} is equivalent to each of the following\/$:$
 \begin{align}
\iim(\pi_{\ell}^{\p} \cdots \pi_1^{\p}) &= \aal_{\ell}^{\p}
	\qquad (\ell=1,\ldots, i)\,, \label{perps-surj}
\\
\iim(\pi_1 \cdots \pi_{\ell}) &= \aal_1
 \qquad (\ell=1,\ldots, i)\,. \label{noperps-surj}
\end{align}
\end{prop}

\begin{proof}
By simple linear algebra, \eqref{cover} is equivalent to each of the
following, where $\vvec{0}$ denotes the zero bundle:
\begin{align}
\aal_{\ell}^{\p} \cap \aal_{\ell-1} = \vvec{0}
	\qquad &(\ell=2,\ldots, i)\,, \label{nondeg-int} \\
\pi_{\ell}^{\p}\aal_{\ell-1}^{\p} = \aal_{\ell}^{\p}
	\qquad &(\ell=2,\ldots, i)\,. \label{nondeg-perps}
\end{align}

Noting that $\iim(\pi_{\ell}^{\p} \cdots \pi_1^{\p}) = \pi_{\ell}^{\p}\!\bigl(\iim(\pi_{\ell-1}^{\p} \cdots \pi_1^{\p})\bigr)$, the
equivalence of \eqref{perps-surj} and \eqref{nondeg-perps} follows by induction; similarly for the equivalence of \eqref{noperps-surj} and \eqref{cover}.
\end{proof}

\begin{prop} \label{prop:cover-props}
Let $\phi_i$ be given by \eqref{phi_i} for some constant map $\phi_0$ and sequence $\aal_1, \ldots, \aal_i$ of unitons which satisfy the covering condition \eqref{cover}.  Then,

{\rm (i)} each $\aal_{\ell}$ is antibasic for $\phi_{\ell-1}$, i.e.,
\begin{equation}
\iim A^{\phi_{\ell-1}}_z \subseteq \aal_{\ell}
	 \qquad (\ell=1,\ldots, i)\,; \label{antibasic}
\end{equation}

{\rm (ii)} we have
$$
\iim A_z^{\phi_{\ell-1}} = \iim A_z^{\phi_{\ell-1}}|_{\aal_{\ell-1}}
= \iim A_z^{\phi_{\ell-1}}|_{\aal_{\ell}}
	\qquad (\ell=2,\ldots, i)\,.
$$
\end{prop}

\begin{proof}
(i) To establish \eqref{antibasic}, we prove the equivalent statement that $\aal_{\ell}^{\p}$ is basic for $\phi_{\ell}$\,, i.e.,
\begin{equation} \label{basic}
\aal_{\ell}^{\p} \subseteq \kker A^{\phi_{\ell}}_z
	\qquad (\ell=1,\ldots, i)\,.
\end{equation}	
To show this, apply \eqref{commute}(i) $\ell$ times giving $A^{\phi_{\ell}}_z \pi_{\ell}^{\p} \cdots \pi_1^{\p} = \pi_{\ell}^{\p} \cdots \pi_1^{ \p} A^{\phi_0}_z$ $=0$.  Together with \eqref{nondeg-perps} this shows that $A^{\phi_{\ell}}_z$ vanishes on $\aal_{\ell}^{\p}$, which establishes \eqref{basic}.

(ii) The first equality holds because $\aal_{\ell-1}^{\p}$ is in the kernel of $A_z^{\phi_{\ell-1}}$, and the second follows from \eqref{cover}.
\end{proof}

The covering condition also has the following vital consequence.

\begin{lem} \label{lem:mero-perp}
With hypothesis as in Proposition \ref{prop:cover-props},
the map $\si \mapsto C^i_i\,\si = \pi_i^{\p}\cdots \pi_1^{\p}(\si)$ maps the set of meromorphic $\CC^n$-valued functions onto the set of meromorphic sections of\/ $(\aal_i^{\p}, D^{\phi_i}_{\zbar})$.
\end{lem}

\begin{proof}  Since $\aal_i^{\p}$ is a uniton for
$\phi_i$, it is a holomorphic subbundle of
$(\CC^n, D^{\phi_i}_{\zbar})$. That $\pi_i^{\p}\cdots \pi_1^{\p}(\si)$ is a meromorphic section of that subbundle for any meromorphic $\si$ follows from \eqref{commute}(ii); that this gives all meromorphic sections follows from \eqref{perps-surj}.
\end{proof}

We now show that our explicit formulae give all harmonic maps obtained from unitons satisfying the covering condition
\eqref{cover}.  In the next section, we shall see that these are
\emph{all} harmonic maps $M^2 \to \U(n)$ of finite uniton number, and thus all harmonic maps from the $2$-sphere to $\U(n)$.
The following result is the key to our explicit formulae.

\begin{prop} \label{prop:calcs}
Let $\phi_i$ be given by \eqref{phi_i} for some integer $i \geq 1$, constant map $\phi_0$, and sequence of unitons $\aal_1,\ldots, \aal_i$.
For any meromorphic vectors $H_0,H_1,\ldots, H_i$\,, set
\begin{equation} \label{K-abstr}
K_{\ell}^{(k)} = \sum_{s=k}^{\ell} C^{\ell}_s H^{(k)}_{s-k}
\qquad (0 \leq k \leq \ell \leq i)\,.
\end{equation}
Suppose that each $K_{\ell}^{(k)}$ is a section of\/
$\aal_{\ell+1}$ for $\ell = 0,\ldots, i-1$. Then

{\rm (i)} $K_i^{(k)}$ is a meromorphic section of\/ $(\CC^n, D^{\phi_i}_{\zbar})\,;$

{\rm (ii)} $A^{\phi_i}_z(K_i^{(k)})
	= \left\{\begin{array}{cl}
		-K_i^{(k+1)}\,, & \mbox{ if \ $k < i$\,$,$} \\
			0\,, 		&  \mbox{ if \ $k=i$\,$;$}
	\end{array}\right.$

{\rm (iii)} $\pi_i(K_i^{(k)}) = K_{i-1}^{(k)}\,.$
\end{prop}

\begin{proof}
(iii) {}From \eqref{C-Pascal} we obtain
$$
K_i^{(k)} = \sum_{s=k}^{i} C^{i}_s H_{s-k}^{(k)}
	= \sum_{s=k}^{i-1} C^{i-1}_s H_{s-k}^{(k)}
		+ \pi_i^{\p} \sum_{s=k}^{i-1} C^{i-1}_{s-1} H_{s-k}^{(k)}\;.
$$
The result follows by noting that the first term is $K_{i-1}^{(k)}$\,, which is in $\aal_i$\,.

\smallskip

Parts (i) and (ii) are proved by some fairly long calculations which we give at the end of this section.
\end{proof}

\begin{thm} \label{thm:all-covering}
Let $\phi_0:M^2 \to \U(n)$ be a constant map, let $r$ be a positive integer, and let $\aal_1, \ldots, \aal_r$ be subbundles of\/ $\CC^n$.  For $i = 1,\ldots, r$, define $\phi_i$ by \eqref{phi_i}.  Then
$\aal_1,\ldots,\aal_r$ is a sequence of unitons satisfying the covering condition \eqref{cover} if and only if they are given by formula \eqref{alpha-H} for some array $(H_{i,j})_{0 \leq i \leq r-1, \ 1\leq j \leq n}$ of meromorphic $\C^n$-valued functions or, equivalently, by formula \eqref{alpha-B} for some holomorphic subbundle $\BB$ of\/ $\CC^{rn}$.
\end{thm}

\begin{proof}
(i) Given a holomorphic subbundle $\BB$ of $\CC^{rn}$,
we define the $\aal_i$ by \eqref{alpha-B}.  {}From
Proposition \ref{prop:calcs} we see inductively that
(i) each $\aal_i^{(k)}$, and so $\aal_i$\,, is a holomorphic subbundle of $(\CC^n,D^{\phi_{i-1}}_{\zbar})$;
(ii) $A^{\phi_{i-1}}_z$ maps $\aal_i^{(k)}$ onto $\aal_i^{(k+1)}$ for $k < i-1$, or to the zero bundle for $k = i-1$, in particular $\aal_i$ is closed under $A^{\phi_{i-1}}_z$.  Hence the $\aal_i$ are unitons; Proposition
\ref{prop:calcs}(iii) show that they satisfy the covering condition \eqref{cover}.
				
(ii) Conversely, suppose that we are given unitons $\aal_i$ satisfying \eqref{cover}; we shall show that they are given
by \eqref{alpha-H} for some array $(H_{i,j})$ of meromorphic vectors.
We shall, in fact, obtain an array in the echelon form of Remark
 \ref{rem:constr-H}(iv) with $\aal_i = \aal_i^{(0)}$ for each $i$.

First, since $\phi_0$ is a constant map, it induces the standard connection $D^{\phi_0}$ on $\CC^n$ so that the uniton $\aal_1$ of $\phi_0$  is a holomorphic subbundle of
$(\CC^n,D^{\phi_0}_{\zbar}) = (\CC^n,\pa_{\zbar})$.
Choose a basis $\{H_{0,1},\ldots, H_{0,d_1}\}$ for it, and set $H_{0,j}$
equal to zero for $j > d_1$. We have thus obtained the first row of our data.  This will not be changed.

Next, $\phi_1$ is formed by adding the uniton $\aal_1$ to
$\phi_0$.  Any uniton $\aal_2$ for $\phi_1$ is a holomorphic
 subbundle of $(\CC^n, D^{\phi_1}_{\zbar})$; to describe
 it, we need to find all meromorphic sections for
$(\CC^n, D^{\phi_1}_{\zbar})$.  {}From \eqref{uniton-DA}(ii), we see that any meromorphic section of $(\CC^n, D^{\phi_0}_{\zbar})$ which lies in $\aal_1$ is also a meromorphic section of $(\CC^n, D^{\phi_1}_{\zbar})$ --- a special property of this first step. On the other hand, by \eqref{commute}(ii), we have meromorphic sections $\pi_1^{\p}H_1$ of $\aal_1^{\p}$ where $H_1$ is an arbitrary meromorphic vector; by Lemma \ref{lem:mero-perp}, all meromorphic sections of $\aal_1^{\p}$ are of this form.  Since $\aal_1^{\p}$ is a holomorphic subbundle of $(\CC^n,D^{\phi_1}_{\zbar})$, these sections are also meromorphic sections of $(\CC^n, D^{\phi_1}_{\zbar})$.  Thus the most general meromorphic section of $(\CC^n, D^{\phi_1}_{\zbar})$ is of the form
$H_0 + \pi_1^{\p}H_1$ where $H_0$ and $H_1$ are meromorphic vectors and $H_0$ is a section of $\aal_1$.

Now we must find a meromorphic basis for $\aal_2$.
By the covering property $\pi_1(\aal_2) = \aal_1$, since
$\{H_{0,j}: 1 \leq j \leq d_1\}$ is a basis for $\aal_1$, for each
$j=1,\ldots, d_1$, there is a meromorphic vector $H_{1,j}$ such that
$H_{0,j} + \pi_1^{\p}H_{1,j}$ is a meromorphic section of $\aal_2$; these are clearly linearly independent.
Additionally, since $\aal_1^{\p} \cap \aal_2$ is a holomorphic subbundle of $\aal_1^{\p}$, we may choose
$H_{1,j}$  $(d_1+1 \leq j \leq d_2)$ with some $d_2 \geq d_1$ such that $\{\pi_1^{\p}H_{1,j}\}$ is a basis for it. Then, recalling that we put
$H_{0,j}$ equal to zero for $j > d_1$, the set
$\{H_{0,j} + \pi_1^{\p}H_{1,j}: 1 \leq j \leq d_2\}$
is a meromorphic basis for $\aal_2$\,.

As before, set $H_{1,j} = 0$ for $j > d_2$; then we have defined the first two rows of $(H_{i,j})$; these will not be changed.

\smallskip

The induction step is slightly more complicated.
Suppose that we have defined the first $i$ rows of $(H_{i,j})$ such that the unitons $\aal_1,\ldots, \aal_i$ are given by \eqref{K} for $j \leq d_i$\,.
We assume that $(H_{i,j})$ is in echelon form; in particular, there are
$d_1 \leq d_2 \leq \cdots \leq d_i$ such that
$H_{\ell,j}$ is zero when $\ell < i$ and $j > d_{\ell+1}$.
The $K_{i-1,j}$ thus give a basis for $\aal_i$ which is meromorphic with respect to
$D^{\phi_{i-1}}_{\zbar}$;
unfortunately, when $i > 1$, this is not usually meromorphic with respect to $D^{\phi_{i}}_{\zbar}$.  However, set
$$
\textstyle \wt{K}_{i,j} = \sum_{s=0}^{i-1}C_s^i H_{s,j}\,;
$$
thus $\wt{K}_{i,j}$ is given by the formula \eqref{K} with $H_{i,j} = 0$.
Then $\wt{K}_{i,j}$ \emph{is}
 meromorphic with respect to $D^{\phi_{i}}_{\zbar}$; further, from
Proposition \ref{prop:calcs},
$\pi_{i-1}(\wt{K}_{i,j}) = K_{i-1,j}$\,.

Next, by Lemma \ref{lem:mero-perp}, any meromorphic section of $\aal_i^{\p}$ is of the form $\pi_i^{\p} \cdots \pi_1^{\p} H_i = C^i_i H_i$ for some meromorphic vector $H_i$\,.

Now, by the covering property $\pi_i(\aal_{i+1}) = \aal_i$\,, since $\{K_{i-1,j} : 1 \leq j \leq d_i\}$ is a basis for $\aal_i$, for each
$j=1,\ldots, d_i$ there must be a meromorphic section of $\aal_{i+1}$ of the form
$K_{i,j} = \wt{K}_{i,j} + C^i_i H_{i,j}$ for some meromorphic vector $H_{i,j}$\,.

Additionally, since $\aal_i^{\p} \cap \aal_{i+1}$ is a holomorphic subbundle of $\aal_i^{\p}$, we can choose meromorphic vectors
$H_{i,j}$ \ $(j=d_i+1, \ldots, d_{i+1})$ for some $d_{i+1} \geq d_i$ such that
$\{\pi_i^{\p} \cdots \pi_1^{\p} H_{i,j} : d_i+1 \leq j \leq d_{i+1}\}$ is a
meromorphic basis for  $\aal_i^{\p} \cap \aal_{i+1}$\,.

Then, recalling that $H_{\ell,j}$ is zero if $\ell < i$ and $j > d_i$\,, the set $\{K_{i,j} : 1 \leq j \leq d_{i+1}\}$ is a meromorphic basis for
$\aal_{i+1}$\,.  Set $H_{i,j}$ equal to zero for $j > d_{i+1}$ and we have defined the first $i+1$ rows of $(H_{i,j})$, completing the
induction step.
\end{proof}

\begin{rem} \label{rem:gen-uniton}
The generating subbundle $\aal_i^{(0)}$ depends on the data $\BB$ (or $(H_{i,j})$); however,
since $\aal_i = \aal_i^{(0)} + \iim A_z^{\phi_{i-1}}$,
it is always a holomorphic subbundle of $\aal_i$ which contains a complement of $\iim A_z^{\phi_{i-1}}$ in $\aal_i$.  The largest such subbundle is $\aal_i$ itself; the data constructed in part (ii) corresponds to that case.

Once $\aal_i^{(0)}$ is defined, the subbundles $\aal_i^{(k)}$ defined
in Remark \ref{rem:constr-H} or \eqref{alpha-B} are characterized by
$\aal_i^{(k)} = \iim A_z^{\phi_{i-1}} \big\vert_{\aal_i^{(k-1)}}$ \ $(k=1,\ldots, i-1)$. In particular, they are independent of the choice of local complex coordinates used to form the derivatives in formulae \eqref{alpha-H} and \eqref{alpha-B}.
\end{rem}

B. \emph{Calculations}

\smallskip

We turn to the proof of Proposition \ref{prop:calcs}(i).
We need the following local calculations which we do in a coordinate domain $(U,z)$ under the hypotheses of Proposition \ref{prop:calcs}.
The reader not interested in the details is invited to skip to Section \ref{sec:fact}.

\begin{lem} \label{lem:Dzbar}
For any holomorphic map $H:U \to \CC^n$, and integers $\ell$, $s$ with $1 \leq \ell \leq i$ and\/ $0\leq s\leq \ell-1$,
$$
D_{\zbar}^{\phi_{\ell}}(\pi_{\ell}^{\p} C_s^{\ell-1}H)=-\pi_{\ell}^{\p}\pa_{\zbar}(C^{\ell-1}_{s+1}H)\,.
$$
\end{lem}

\begin{proof}
For $s=\ell-1$, the right-hand side is zero and the lemma follows from \eqref{commute}(ii).  We now assume that $s \leq \ell-2$ and proceed by induction on $\ell$. The first case is
$\ell=2$ and $s=0$; in that case, using \eqref{commute}(ii) and holomorphicity of $H$, we have
$$
D_{\zbar}^{\phi_2}\pi_2^{\p} C_0^1 H
	=\pi_2^{\p} D_{\zbar}^{\phi_1}H
=\pi_2^{\p}(D_{\zbar}^{\phi_0}H-\pa_{\zbar}\pi_1^{\p}H 		+\pi_1^{\p}\pa_{\zbar}H)
=-\pi_2^{\p} \pa_{\zbar}(C_1^1 H).
$$
Next, suppose that the result holds up to some value of $\ell$ with $2 \leq \ell \leq i-1$; we shall show that it holds for $\ell+1$.
Using \eqref{commute}(ii), \eqref{C-Pascal} and the induction hypothesis,
\begin{align*}
D_{\zbar}^{\phi_{\ell+1}}\pi_{\ell+1}^{\p} C_{s}^{\ell}H&=\pi_{\ell+1}^{\p} D_{\zbar}^{\phi_{\ell}}C_s^{\ell}H=\pi_{\ell+1}^{\p} (D_{\zbar}^{\phi_{\ell}}\pi_{\ell}^{\p} C_{s-1}^{\ell-1}H+D_{\zbar}^{\phi_{\ell}}C_s^{\ell-1}H)
\\
&=\pi_{\ell+1}^{\p}(-\pi_{\ell}^{\p}\pa_{\zbar}C_s^{\ell-1}H+D_{\zbar}^{\phi_{\ell-1}}C_s^{\ell-1}H-\pa_{\zbar}\pi_{\ell}^{\p} C_s^{\ell-1}H+\pi_{\ell}^{\p}\pa_{\zbar}C_s^{\ell-1}H)
\\
&=\pi_{\ell+1}^{\p} D_{\zbar}^{\phi_{\ell-1}}C_s^{\ell-1}H-\pi_{\ell+1}^{\p}\pa_{\zbar}\pi_{\ell}^{\p} C_s^{\ell-1}H
\intertext{(repeating the argument with the first term)}
&=\pi_{\ell+1}^{\p} (D_{\zbar}^{\phi_{\ell-1}}\pi_{\ell-1}^{\p} C_{s-1}^{\ell-2}H+D_{\zbar}^{\phi_{\ell-1}} C_{s}^{\ell-2}H)-\pi_{\ell+1}^{\p}\pa_{\zbar}\pi_{\ell}^{\p} C_s^{\ell-1}H
\\
&=\pi_{\ell+1}^{\p}(-\pi_{\ell-1}^{\p} \pa_{\zbar}C_{s}^{\ell-2}H
+D_{\zbar}^{\phi_{\ell-2}}C_s^{\ell-2}H-\pa_{\zbar}\pi_{\ell-1}^{\p} C_s^{\ell-2}H\\
&\qquad +\pi_{\ell-1}^{\p}\pa_{\zbar}C_s^{\ell-2}H)
- \pi_{\ell+1}^{\p}\pa_{\zbar}\pi_{\ell}^{\p} C_s^{\ell-1}H
\\
&=\pi_{\ell+1}^{\p} D_{\zbar}^{\phi_{\ell-2}}C_s^{\ell-2}H-\pi_{\ell+1}^{\p}\pa_{\zbar}(\pi_{\ell-1}^{\p} C_s^{\ell-2}H+\pi_{\ell}^{\p} C_s^{\ell-1}H)
\\
\intertext{(repeating the argument another $\ell-s-2$ times)}
&=\pi_{\ell+1}^{\p} D_{\zbar}^{\phi_s}C_s^s H-
\\
-\pi_{\ell+1}^{\p}&\pa_{\zbar}(\pi_{s+1}^{\p} C_s^s H +\pi_{s+2}^{\p} C_s^{s+1}H+...+\pi_{\ell-1}^{\p} C_s^{\ell-2}H+\pi_{\ell}^{\p} C_s^{\ell-1}H).
\end{align*}

{}From Lemma \ref{lem:mero-perp}, the first term in the above expression vanishes. Using \eqref{C-Pascal} to add up the remaining terms, we obtain
$$
D_{\zbar}^{\phi_{\ell+1}}\pi_{\ell+1}^{\p} C_{s}^{\ell}H=
-\pi_{\ell+1}^{\p}\pa_{\zbar}(C_{s+1}^{\ell}H),
$$
completing the proof.
\end{proof}

\emph{Proof of Proposition \ref{prop:calcs}(i).}
We prove part (i) of Proposition \ref{prop:calcs} for $k=0$.
  For other values of $k$, it will follow from that case and the holomorphicity of $A^{\phi_{i-1}}_z$.

It clearly holds for $i=0$.  Suppose that it holds for some value of
$i \geq 0$.  Using \eqref{C-Pascal} and
$K_{\ell}^{(0)} \in \Ga(\aal_{\ell+1})$, we have\\[-3ex]
\begin{align*}
D^{\phi_{i+1}}_{\zbar}K_{i+1}^{(0)}&=D^{\phi_{i+1}}_{\zbar}(\sum_{s=0}^{i+1}C_s^{i+1}H_{s})\\
&=D^{\phi_{i+1}}_{\zbar}\Bigl(K_i^{(0)}+\pi_{i+1}^{\p} \sum_{s=0}^{i}C^i_{s} H_{s+1} \Bigr)\\[-1ex]
\intertext{(using \eqref{uniton-DA}(ii), Lemma \ref{lem:Dzbar}) and the induction hypothesis)}\\[-5ex]
&=D^{\phi_{i}}_{\zbar}K_i^{(0)} -\pa_{\zbar}(\pi_{i+1}^{\p})K_i^{(0)}-\pi_{i+1}^{\p}\pa_{\zbar} \sum_{s=0}^{i-1} C^i_{s+1} H_{s+1}
\\[-2ex]
&=0+\pi_{i+1}^{\p}\pa_{\zbar}
	\bigl(K_i^{(0)}-\sum_{s=1}^{i} C^i_{s} H_{s}\bigr)
\\[-1ex]
&=\pi_{i+1}^{\p}(\pa_{\zbar}H_{0})=0\,,
\end{align*}\\[-2.5ex]
completing the induction step.
\qed

\medskip

To prove Proposition \ref{prop:calcs}(ii) we need a general algebraic lemma.
\begin{lem} \label{lem:Az-abstract}
For any integers $\ell$, $k$ with  $0 \leq k \leq \ell \leq i$ and smooth functions $J_0,J_1, \ldots, J_{\ell-k}$\,,
$$
A^{\phi_{\ell}}_z \Bigl(\sum_{s=k-1}^{\ell} C_s^{\ell} J_{s-k+1} \Bigr) -\pa_z \sum_{s=k}^{\ell}C_s^{\ell}J_{s-k} = - \sum_{s=k}^{\ell}C_s^{\ell}J^{(1)}_{s-k}\,.
$$
\end{lem}

\begin{proof}
For $\ell=0$, this is trivially true. Suppose that it is true for some value of $\ell$ with $0 \leq \ell \leq i-1$.  Then, using \eqref{C-Pascal}, \eqref{commute}(i) and \eqref{uniton-DA}(i),
\begin{align*}
A^{\phi_{\ell+1}}_z &\Bigl(\sum_{s=k-1}^{\ell+1} C_s^{\ell+1} J_{s-k+1} \Bigr) -\pa_z \sum_{s=k}^{\ell+1} C_s^{\ell+1} J_{s-k}
\\
	=& \ A^{\phi_{\ell+1}}_z \Bigl(\sum_{s=k-1}^{\ell+1}
(\pi_{\ell+1}^{\p} C_{s-1}^{\ell} + C_s^{\ell}) J_{s-k+1} \Bigr)
-\pa_z \sum_{s=k}^{\ell+1}(\pi_{\ell+1}^{\p} C_{s-1}^{\ell}+C_s^{\ell})J_{s-k}
\\
=& \ \pi_{\ell+1}^{\p}A_z^{\phi_{\ell}} \Bigl(\sum_{s=k-2}^{\ell} C_s^{\ell} J_{s-k+2} \Bigr) + A_z^{\phi_{\ell}}
	\Bigl(\sum_{s=k-1}^{\ell} C_s^{\ell} J_{s-k+1} \Bigr)
\\
&+ (\pa_z\pi_{\ell+1}^{\p})
	\Bigl(\sum_{s=k-1}^{\ell} C_s^{\ell} J_{s-k+1} \Bigr)
-\pa_z\Bigl(\pi_{\ell+1}^{\p}\!\!\sum_{s=k-1}^{\ell} C_s^{\ell} J_{s-k+1} \Bigr)
	-\pa_z \Bigl(\sum_{s=k}^{\ell} C_s^{\ell} J_{s-k} \Bigr)
\\
=& \ \pi_{\ell+1}^{\p}\Bigl(A_z^{\phi_{\ell}}\!\!\!
	\sum_{s=k-2}^{\ell} C_s^{\ell} J_{s-k+2}
	-\pa_z \!\!\sum_{s=k-1}^{\ell}C_s^{\ell}J_{s-k+1}\Bigr)\\
&\qquad +	A_z^{\phi_{\ell}}\!\!\!
	\sum_{s=k-1}^{\ell} C_s^{\ell} J_{s-k+1}
	-\pa_z \!\sum_{s=k}^{\ell}C_s^{\ell}J_{s-k}
\\
\intertext{(using the induction hypothesis twice)}
=& \ -\pi_{\ell+1}^{\p}\sum_{s=k-1}^{\ell}C_s^{\ell}J^{(1)}_{s-k+1} - \sum_{s=k}^{\ell}C_s^{\ell}J^{(1)}_{s-k}\\
=& \ -\sum_{s=k}^{\ell+1} (\pi_{\ell+1}^{\p}C_{s-1}^{\ell} + C_s^{\ell})J^{(1)}_{s-k}
	= -\sum_{s=k}^{\ell+1} C_s^{\ell+1} J^{(1)}_{s-k}\,,
\end{align*}
which establishes the induction step.
\end{proof}

\emph{Proof of Proposition \ref{prop:calcs}(ii).}
Define $K_i^{(k)}$ to be zero if $k > i$, then we can treat the cases $i < k$ and $i = k$ together.
For $i=0$, the result is trivial.  Suppose that it is true for some value of $i \geq 0$.  Then, using \eqref{C-Pascal}, \eqref{commute}(i) and \eqref{uniton-DA}(i),
\begin{align*}
A_z^{\phi_{i+1}}& K_{i+1}^{(k)}
= A_z^{\phi_{i+1}} \Bigl(\pi_{i+1}^{\p} \sum_{s=k}^{i+1} C_{s-1}^{i} H_{s-k}^{(k)} + \sum_{s=k}^{i} C_s^{i} H_{s-k}^{(k)} \Bigr)
\\
=& \ \pi_{i+1}^{\p} A_z^{\phi_{i}} \sum_{s=k}^{i+1} C_{s-1}^{i} H_{s-k}^{(k)}
+ A_z^{\phi_{i}} K_{i}^{(k)}
+ \pa_z(\pi_{i+1}^{\p}K_{i}^{(k)})
	- \pi_{i+1}^{\p}(\pa_z K_{i}^{(k)})\\
\intertext{(using $\pi_{i+1}^{\p}K_{i}^{(k)} = 0$, \eqref{C-Pascal} and the induction hypothesis)}
=& \ \pi_{i+1}^{\p}\Bigl(A_z^{\phi_{i}}\!\! \sum_{s=k-1}^{i} C_s^{i}H_{s-k+1}^{(k)} - \pa_z \sum_{s=k}^{i} C_s^{i} H_{s-k}^{(k)} \Bigr) - K_{i}^{(k+1)}\\
\intertext{(using Lemma \ref{lem:Az-abstract} and then \eqref{C-Pascal})\,}
=& \ -\pi_{i+1}^{\p}\Bigl(\sum_{s=k}^{i} C_s^{i}H_{s-k}^{(k+1)} \Bigr) - K_{i}^{(k+1)}\\
=& - \sum_{s=k+1}^{i+1} (\pi_{i+1}^{\p}C_{s-1}^{i}+C_s^{i})H_{s-k-1}^{(k+1)}
 = -\sum_{s=k+1}^{i+1} C_s^{i+1} H_{s-k-1}^{(k+1)} = -K_{i+1}^{(k+1)},
\end{align*}
completing the induction step.
\qed

\section{Factorization of an arbitrary harmonic map}\label{sec:fact}

In this section we show that any harmonic map of finite uniton number from $M^2$ to $\U(n)$ can be obtained as the product \eqref{fact} of unitons satisfying the covering condition \eqref{cover}, and is thus given explicitly by our formulae.   For this we need to recall Uhlenbeck's concept of \emph{extended solution} \cite{Uhl}.

As in the last section, for any smooth map $\phi:M^2 \to \U(n)$ we define a $\uu(n)$-valued $1$-form $A = A^{\phi}$ by
\begin{equation} \label{A}
A = \half \phi^{-1}\dd\phi\,;
\end{equation}
this is one half the pull-back of the Maurer--Cartan form.  Clearly, $A$ satisfies the pull-back of the Maurer--Cartan equations:
\begin{equation} \label{M-C}
\dd A + [A,A] = 0\,,
\end{equation}
where $[A,A]$ is the $\uu(n)$-valued $2$-form defined at each point $p \in M$ by
$[A,A](X,Y) = [A(X),A(Y)]$ \ $(X,Y \in T_pM)$.

Conversely, given a smooth $\uu(n)$-valued $1$-form $A$, there exist locally (globally, if $M^2$ is simply connected) smooth maps $\phi:M^2 \to \U(n)$ satisfying \eqref{A} if and only if \eqref{M-C} holds, i.e., \eqref{M-C} is the \emph{integrability condition} for the equation \eqref{A}.

Now, on writing $A = A_z \dd z + A_{\zbar} \dd\zbar$, \eqref{M-C} read
\begin{equation} \label{M-C-cx}
\pa_{\zbar}A_z - \pa_z A_{\zbar} + 2[A_{\zbar}, A_z] = 0\,.
\end{equation}
Further, the harmonic equation \eqref{hol} for $\phi$ can be written as either of the equivalent equations:
\begin{equation} \label{ha-cx}
\textrm{(a) } \ \pa_{\zbar}A_z + [A_{\zbar}, A_z] = 0\,,
\quad
\textrm{(b) } \ \pa_z A_{\zbar} + [A_z, A_{\zbar}] = 0\,.
\end{equation}

For each $\la \in S^1$, set
\begin{equation} \label{loop}
A_{\la} = \half(1-\la^{-1})A_z \dd z + \half(1-\la)A_{\zbar}\dd\zbar\,;
\end{equation}
then on calculating the integrability condition $\dd A_{\la} + [A_{\la}, A_{\la}] = 0$ and equating coefficients of $\la$, we see that $A_{\la}$ is integrable for each $\la$ if and only if equations \eqref{M-C-cx} and \eqref{ha-cx} all hold, i.e, if and only if $A = \half \phi^{-1}\dd\phi$ for some harmonic map $\phi:M^2 \to \U(n)$.
In that case we can find, at least locally, an $S^1$-family of smooth maps
$\Phi = \Phi_{\la}:M^2 \to \U(n)$
with
\begin{equation} \label{M-C-E}
\half \Phi_{\la}^{-1} \dd\Phi_{\la} = A_{\la}
\qquad (\la \in S^1)\,,
\end{equation}
equivalently,
\begin{equation} \label{M-C-E2}
\pa_z\Phi_{\la} = (1-\la^{-1})\Phi_{\la} A^{\phi}_z
\quad \text{and} \quad
\pa_{\zbar}\Phi_{\la} = (1-\la)\Phi_{\la} A^{\phi}_{\zbar}\,.
\end{equation}
We are thus led to the following definition.
\begin{defn} \cite{Uhl}
Say that $\Phi = \Phi_{\la}:M^2 \to \U(n)$ is an \emph{extended solution} (for $\phi$) if it satisfies \eqref{M-C-E} for each
$\la \in S^1$, with $A_{\la}$ given by \eqref{loop} for some harmonic map $\phi:M^2 \to \U(n)$ and
\begin{equation} \label{Phi1}
\Phi_1(p) = I \text{ for all } p \in M^2,
\end{equation}
where $I$ is the identity of the group.
\end{defn}

Note that our extended solutions will be globally defined on $M^2$.  {}From equations \eqref{loop} and \eqref{M-C-E}, $\Phi_1$
is necessarily constant so that the condition \eqref{Phi1} can be achieved; it means that $\Phi$ can be interpreted as a map into a loop group, see Section \ref{sec:loop}.
Note that any two extended solutions for a harmonic map differ by a function (`constant loop') $Q:S^1 \to \U(n)$ with $Q(1) = 1$.
Further, $\Phi_{-1}$ is left-equivalent to $\phi$, i.e.,
$\Phi_{-1} = A \phi$ for some constant $A \in \U(n)$; we do not insist that $\Phi_{-1}$ equal $\phi$\,: when coupled with the condition \eqref{Phi1}, this is not always convenient,
as the following example shows.

\begin{exm} \label{ex:extended-soln}
An extended solution for the harmonic map $\phi:M^2 \to \U(n)$
 given by \eqref{fact} is
\begin{equation} \label{fact-ext}
\Phi = (\pi_1 + \la\pi_1^{\p})
	\cdots (\pi_r + \la \pi_r^{\p})\,.
\end{equation}
We see that $\phi = \phi_0 \Phi_{-1}$.   An extended
solution with $\Phi_{-1} = \phi$ would need
the extra factor
$\half(I+\phi_0) + \half\la(I-\phi_0)$.
\end{exm}

Let $\gl(n,\C)$ denote the Lie algebra of
$n \times n$ matrices; this is
the complexification of $\uu(n)$.
The extended solution extends to a family of maps
$\Phi_{\la}:M^2 \to \gl(n,\C)$
with $\Phi_{\la}$ a holomorphic function of
$\la \in \C \setminus \{0\}$.
Hence it can be expanded as a Laurent series,
$\Phi = \sum_{i=-\infty}^{\infty} \la^i T_i$\,, where each $T_i = T_i^{\Phi}$ is a smooth map from
$M^2$ to $\gl(n,\C)$.  That $\Phi_{-1}$, and so
$\Phi_{\la}$ $(\la \in S^1)$, has values in $\U(n)$ is expressed by the
\emph{reality condition}
\begin{equation} \label{reality0}
\Phi_{\la}\Phi_{\la}^* = \Phi_{\la}^*\Phi_{\la} = I
	\qquad (\la \in S^1) \,.
\end{equation}

Recall that a harmonic map $\phi:M^2 \to \U(n)$ is said to be
\emph{of finite uniton number} if it can be written as a
 product \eqref{fact} of a finite number of unitons; equivalently \cite{Uhl}, it has a \emph{polynomial} extended solution:
\begin{equation} \label{ext-sol}
\Phi = T_0 + \la T_1 + \cdots + \la^i T_i\,.
\end{equation}
The \emph{(minimal) uniton number} of $\phi$ is the least number of unitons such that $\phi$ can be written as a product \eqref{fact};
equivalently, it is the least degree of all its polynomial extended solutions.
Note that, if $T_0$ is zero, the polynomial is divisible by the constant loop $\la$, so an extended solution of least degree $r$ will have both $T_0$ and $T_r$ non-zero.

Note also that any extended solution which is a Laurent polynomial $\sum_{\ell=-s}^{\ell=t} \la^{\ell} T_{\ell}$ can be converted to a polynomial extended solution by multiplication by a constant loop $\la^s$.

Given any polynomial extended solution \eqref{ext-sol}, on equating coefficients, we see that the extended solution equations \eqref{M-C-E2} are equivalent to
\begin{equation} \label{M-C-coeffs}
{\rm (a)} \ \ \pa_z T_{\ell} = (T_{\ell} -  T_{\ell+1}) A^{\phi}_z
\quad \text{and} \quad
{\rm (b)} \ \
\pa_{\zbar} T_{\ell} = (T_{\ell} - T_{\ell-1}) A^{\phi}_{\zbar}
\quad
\forall\, \ell,
\end{equation}
where we set $T_{\ell}$ equal to zero for $\ell<0$ and $\ell>i$.

By equating coefficients in \eqref{reality0} we obtain \emph{reality conditions} on the top and bottom coefficients $T_0$ and $T_i$ in
\eqref{ext-sol}:
\begin{equation} \label{reality}
\textrm{(a)} \ \ T_0 T_i^{\:*} = 0\,, \quad
	\textrm{(b)} \ \ T_i^{\:*} T_0 = 0\,.
\end{equation}

{}From the extended solution we can construct two unitons as follows.

\begin{prop} \label{prop:T0n} Let $\Phi$ be an extended solution \eqref{ext-sol} for a harmonic map $\phi$.  Then
{\rm (i)}  $\iim T_r^{\: *}$ is a basic uniton for $\phi;$
{\rm (ii)} $\kker T_0$ is an antibasic uniton for $\phi$\/.
\end{prop}

\begin{proof}
(i) As in \cite{Uhl,HeSh}, from \eqref{M-C-coeffs}(a) for $i = r$,
we obtain $\pa_z T_r = T_r A^{\phi}_z$ which is equivalent to $\pa_{\zbar}T_r^{\: *} + A^{\phi}_{\zbar}T_r^{\: *} = 0$. It follows that   $T_r^{\: *}$ is a holomorphic endomorphism from $(\CC^n,\pa_{\zbar})$ to $(\CC^n,D^{\phi}_{\zbar})$, so that its image is a holomorphic subbundle of $(\CC^n, D^{\phi}_{\zbar})$.  Further, from \eqref{M-C-coeffs}(b) with $i=r+1$, we obtain $T_r \circ A^{\phi}_{\zbar} = 0$, equivalently $A^{\phi}_z \circ T_r^{\:*} = 0$, so that $\iim T_r^{\: *}$ is basic.

(ii) Similarly, $T_0$ is a holomorphic endomorphism from $(\CC^n,D^{\phi}_{\zbar})$ to $(\CC^n,\partial_{\zbar})$; the rest
is established by Uhlenbeck \cite[\S 14]{Uhl}.
\end{proof}

We now give a unique factorization theorem for extended solutions into unitons satisfying the covering condition. To achieve uniqueness we insist that each uniton $\aal_i$ be \emph{proper} in the sense that $\aal_i$ is neither the zero subbundle $\vvec{0}$ nor the full bundle $\CC^n$.

\begin{prop} \label{prop:fact-ext}
Let $\Phi$ be an polynomial extended solution \eqref{ext-sol} of some degree $r \geq 1$, with non-zero
constant term $T_0^{\Phi}$.  Then there exists a unique factorization
\eqref{fact-ext} with $\aal_1,\ldots, \aal_r$ proper unitons satisfying the covering condition \eqref{cover}.

Further,  write
$\Phi_i = (\pi_1 + \la\pi_1^{\p})
	\cdots (\pi_i + \la \pi_i^{\p})$ \
$(i = 0,1,\dots, r)$.
Then $\Phi_i$ is a polynomial extended solution\/$:$
\begin{equation} \label{Phi_i}
\Phi_i = T^{\Phi_i}_0 + \la T^{\Phi_i}_1 + \cdots + \la^i T^{\Phi_i}_i
\end{equation}
of degree precisely $i$ with non-zero constant term $T_0^{\Phi_i}$.  In fact,
$$
\iim T^{\Phi_i}_0 = \aal_1 \quad
	\mbox{and} \quad
	\iim(T^{\Phi_i}_i)^* = \aal_i^{\p}.
$$	
\end{prop}

\begin{proof}
We shall define extended solutions \eqref{Phi_i} inductively for $i = r, r-1, r-2, \ldots, 1, 0$.  To start the induction, set $\Phi_r = \Phi$.  Now suppose that, for some $i$, we have defined an extended solution $\Phi_i$ of degree precisely $i$ with non-zero constant term.  Set $\aal_i = \kker T_i^{\Phi_i}$, equivalently $\aal_i^{\p} = \iim(T_i^{\Phi_i})^*$.
By Proposition \ref{prop:T0n},
$\aal_i^{\p}$ is a basic uniton for $\Phi_i$.  Now $T_i^{\Phi_i}$ and $T_0^{\Phi_i}$ are non-zero, and from the reality conditions \eqref{reality}, $T_i$ cannot be an isomorphism, thus $\aal_i^{\p}$ is proper.  Use it to obtain the extended solution:
\begin{equation} \label{inverse-step}
\Phi_{i-1} = \Phi_i(\pi_i + \la^{-1} \pi_i^{\p})
	= \sum_{\ell=-1}^i \la^{\ell} T_{\ell}^{\Phi_{i-1}}\,;
\end{equation}
we shall show that this is polynomial of degree precisely $i-1$ with non-zero constant term.

To do this, note that the coefficient of $\la^{-1}$ in
 \eqref{inverse-step} is $T_0^{\Phi_i} \pi_i^{\p}$;
since $\aal_i^{\p} = \iim(T_i^{\Phi_i})^*$,
this is zero by the reality condition \eqref{reality}(a)
for $\Phi_i$. Similarly the coefficient of $\la^i$ in
\eqref{inverse-step} is $T_i^{\Phi_i} \pi_i$, which is
zero since $\aal_i = \kker T_i^{\Phi_i}$.  Hence $\Phi_{i-1}$ is a polynomial of degree at most $i-1$.  Next note that
\eqref{inverse-step} is equivalent to
\begin{equation} \label{step}
\Phi_i = \Phi_{i-1} (\pi_i + \la\pi_i^{\p})\,.
\end{equation}

{}From this equation, $T_i^{\Phi_i} = T_{i-1}^{\Phi_{i-1}} \pi_i$; since, by the inductive hypothesis, $T_i^{\Phi_i}$ is non-zero, so is $T_{i-1}^{\Phi_{i-1}}$, hence $\Phi_{i-1}$ has degree precisely $i-1$. Similarly,
$T_0^{\Phi_i} = T_0^{\Phi_{i-1}} \pi_i^{\p}$, so that the constant term $T_0^{\Phi_{i-1}}$ is non-zero.

In particular, $\Phi_0$ is a non-zero constant polynomial, which must be equal to $I$ by \eqref{Phi1}, so from \eqref{step} we obtain a factorization \eqref{fact-ext}.  Comparing with \eqref{ext-sol} we see that $(T_i^{\Phi_i})^*=\pi_i^{\p} \cdots \pi_1^{\p}$.  Now, by definition, the image of $(T_i^{\Phi_i})^*$ is equal to $\aal_i^{\p}$.  By Proposition \ref{prop:cover-conds}, this implies the covering condition \eqref{cover}.

Again, from \eqref{Phi_i}, we see from $T_0^{\Phi_i} = \pi_1\cdots \pi_i$\,; the covering condition implies that $T_0^{\Phi_i}$ has image $\aal_1$\,.

As for uniqueness, from Proposition \ref{prop:cover-conds}, we see that, for any factorization into proper unitons $\aal_i$ satisfying the covering condition, equation \eqref{perps-surj} of Proposition \ref{prop:cover-conds} tells us that
$\aal_i^{\p} = \iim (T_i^{\Phi_i})^*$ for all $i$. Since the $\Phi_i$ are uniquely determined by $\Phi$ by the inductive process above, it follows that the $\aal_i$ are also uniquely determined by $\Phi$.
\end{proof}

\begin{rem} \label{rem:dual}
The factorization we use is dual to that of Uhlenbeck \cite[Theorem 14.6]{Uhl}, who sets $\aal_i^{\p}$ equal to $\kker T_0^{\Phi_i}$, an antibasic uniton. Indeed, taking the orthogonal complement of each uniton and changing the complex structure converts the one factorization into the other.  The conditions in these factorizations were simplified in \cite{HeSh} to conditions involving $T_0$ or $T_i$\,.
The covering condition \eqref{cover} is equivalent to their condition, but involves only the unitons.
\end{rem}

To apply this to harmonic maps, we recall a condition in \cite{Uhl} discussed in \cite{HeSh}.  Say that a subbundle $\aal$ of $\CC^n$, or the corresponding map $\al:M^2 \to G_*(\C^n)$
to a Grassmannian, is \emph{full} if $\aal$ is not contained in any proper constant subbundle of $\CC^n$.  Say that a polynomial extended solution \eqref{ext-sol} is \emph{of type one} if $\iim T_0$ is full; note that for a harmonic map $\phi=\phi_0(\pi_1-\pi_1^{\p})$ of uniton number one, this is equivalent to the condition that the holomorphic subbundle $\aal_1$ be full.

Uhlenbeck proves the following.

\begin{prop} \label{prop:type1}
{\rm \cite[Theorems 13.2 and 13.3]{Uhl}}
Given a harmonic map $\phi:M^2 \to \U(n)$ of finite uniton number, there is a unique polynomial extended solution of type one.  Furthermore, this has degree equal to the (minimal) uniton number of $\phi$.
\end{prop}

\emph{Sketch Proof of Existence.}
Let $\Phi$ be a polynomial extended solution of degree
$r = $ minimal uniton number of $\phi$, but not of type one; let $\ul{A}$ be a proper constant
 subbundle containing $\iim T_0^{\Phi}$.
Then another polynomial extended solution of degree $r$ for $\phi$ is given by
$(\pi_A + \la^{-1} \pi_A^{\p})\Phi$. Uhlenbeck shows that
 repeating this process a finite number of times will give a polynomial extended solution of type one.
\qed

\begin{exm} \label{exm:non-full}
Suppose that
$h:M^2 \to G_*(\C^n)$ is holomorphic \emph{but not full};
 then the corresponding subbundle $\hh$ lies fully in some proper constant subbundle $\ul{A}$ of $\CC^n$.  As usual,
$\Phi = \pi_h + \la\pi_h^{\p}$ is an extended solution for
$h$, but it is not of type one.  However,
\emph{$\ul{\wt{h}} = \hh \oplus \ul{A}^{\p}$ is full, and its extended solution $\wt{\Phi} = \pi_{\wt{h}} + \la \pi_{\wt{h}}^{\p}$ is of type one}.

Now $\pi_{\wt{h}} - \pi_{\wt{h}}^{\p}$ and $\pi_h -\pi_h^{\p}$
 are the same up to left multiplication by the constant
 $\pi_A-\pi_A^{\p}$; in fact, as in the above sketch proof,
their extended solutions satisfy
$\wt{\Phi} = (\pi_A + \la^{-1} \pi_A^{\p})\Phi$.
Hence,
\emph{$\pi_{\wt{h}} + \la \pi_{\wt{h}}^{\p}$ is the unique polynomial extended solution of type one for both $h$ and $\wt{h}$}.
\end{exm}

As a consequence of the last proposition, we obtain the following unique factorization theorem for harmonic maps.	

\begin{thm} \label{thm:unique-fact}
Given a harmonic map $\phi:M^2 \to \U(n)$ of (minimal) uniton number $r$, there is a unique factorization \eqref{fact} into proper unitons
 $\aal_1,\ldots,\aal_r$ satisfying the covering condition \eqref{cover} with $\aal_1$ full.

All harmonic maps of finite uniton number, and so \emph{all} harmonic maps from $S^2$, are given this way for some $r \leq n-1$.
\end{thm}

\begin{proof}
Let $\Phi$ be the unique polynomial extended solution of type one for $\phi$; by Proposition \ref{prop:fact-ext}, this has a unique factorization by proper unitons satisfying \eqref{cover}; on putting $\la = -1$ this gives a factorization \eqref{fact}.

Further, since the image of
$T_0^{\Phi_i}$ is the same for all $i$, we see that each $\Phi_i$ is of type one; on putting $i=1$ we see that $\aal_1$ is full.

Conversely, given a factorization \eqref{fact} of a harmonic map $\phi$, we obtain a factorization \eqref{fact-ext} of an extended solution $\Phi$; as just seen, fullness of $\aal_1$ implies that $\Phi$ is of type one showing that the factorization of $\phi$ is unique.
\end{proof}

\emph{Proofs of Theorems \ref{thm:constr-H} and \ref{thm:constr-B}.}
Simply combine the last theorem with Theorem \ref{thm:all-covering},
and recall that Uhlenbeck shows that the minimal uniton
number of any harmonic map $M^2 \to \U(n)$ of finite uniton number is at most $n-1$.

\qed

\begin{cor} \label{cor:uniton-no}
Let $\phi_i$ \ $(i = 0,1,\ldots )$ be a sequence of harmonic maps given by \eqref{phi_i} for some constant map
$\phi_0$ and sequence of unitons $\aal_1,\ldots, \aal_i$
which satisfy the covering condition
\eqref{cover}.  Suppose that $\aal_1$ is full.  Then
$\phi_i$ has uniton number precisely $i$.
\qed \end{cor}

\begin{rem} \label{rem:non-unique}
Without fullness of $\aal_1$, the factorization
\eqref{fact} may not be unique.  For example, let
$\al_1 = h:M^2 \to G_1(\C^n) = \CP^{n-1}$ be a non-full holomorphic map as in Example \ref{exm:non-full},
$H_{0,1}$ any meromorphic section of $\hh = \aal_1$, and $H_{1,1}$ any meromorphic section of $\ul{A}^{\p}$, then the uniton $\aal_2$ determined by this data is spanned by $H_{0,1}+\pi_h^{\p}H_{1,1} = H_{0,1} +H_{1,1}$ and $\pi_h^{\p}(H_{0,1}^{(1)})$.
The resulting harmonic map $\phi_2$ has extended solution
\begin{equation} \label{quadr-ex}
\Phi_2 = (\pi_1+\la\pi_1^{\p})(\pi_2+\la\pi_2^{\p})
	= (\pi_h+\la \pi_h^{\p})(\pi_{\al_2}+\la\pi_{\al_2}^{\p})\,,
\end{equation}
which is not of type one.  Carrying out the process in the proof of Proposition \ref{prop:type1}, we get a new extended solution
$(\pi_A + \la^{-1}\pi_A^{\p})\Phi_2$ which simplifies to
the linear polynomial $\pi_g+\la\pi_g^{\p}$ where $g:M^2 \to G_1(\C^n)$ is the holomorphic map spanned by $H_{0,1}+H_{1,1}$\,.  Hence the harmonic map $\phi_2$ has two factorizations into proper unitons satisfying the covering condition, namely that from \eqref{quadr-ex} into two unitons:
$\phi_2 = (\pi_h-\pi_h^{\p})(\pi_{\al_2}-\pi_{\al_2}^{\p})$
and that into one uniton
$\phi_2 = \phi_0(\pi_g-\pi_g^{\p})$ with
$\phi_0$ equal to the constant $\pi_A -\pi_A^{\p}$.
The latter shows that, up to left-equivalence, $\phi_2$ is actually a holomorphic map into a Grassmannian of uniton number one.
This example also shows that Corollary \ref{cor:uniton-no}
 is false without fullness of $\aal_1$.
\end{rem}

In the following examples,  the classification is up to left equivalence, and maps into a Grassmannian $G_*(\C^n)$ are composed with the Cartan embedding into $\U(n)$.

\begin{exm} \label{exm:low-dim}
(i) As is well known, a harmonic map $\phi:M^2 \to \U(n)$
is of \emph{uniton number zero} if and only if it is constant.
Any harmonic map from $M^2 \to \U(1)$ of finite uniton number is constant.

\medskip

(ii) As is also well known, a harmonic map $\phi:M^2 \to \U(n)$ is of \emph{uniton number one} if and only if it is a holomorphic map into a
 Grassmannian $G_{d_1}(\C^n)$, where $1 \leq d_1 = \rk \aal_1 < n$.  If $\aal_1$ is not full, then, as in
 Example \ref{exm:non-full}, we can replace it by a larger full subbundle.

Any non-constant harmonic map of finite uniton number $M^2 \to \U(2)$ is of uniton number one, and so is a holomorphic map $M^2 \to \CP^1$.

 \medskip

(iii) Let $\phi:M^2 \to \U(3)$ be a non-constant harmonic
 map of finite uniton number.  Then, \emph{either}

(a) it has uniton number one and is given by a holomorphic map
$\phi:M^2 \to G_{d_1}(\C^3)$ where $d_1 = 1$ or $2$; \emph{or}

(b) it has uniton number two and is given by \eqref{fact} with unitons
$\aal_1$, $\aal_2$ have rank one and two respectively, with $\aal_1$ full.
The data of Theorem \ref{thm:constr-H} consists of a single column $(H_{0,1},H_{1,1})^T$.
Then
\begin{equation} \label{rank12}
\aal_1 = \spn\{H_{0,1}\} \quad \text{and} \quad
\aal_2
= \spn\{H_{0,1} + \pi_1^{\p}H_{1,1}\,,\, \pi_1^{\p}H_{0,1}^{(1)}\,\}\,.
\end{equation}
Because of the low dimensions, this can be written as
$$
\aal_1 = \spn\{H_{0,1}\} \quad \text{and} \quad
\aal_2 = \spn\{H_{0,1} + \mu G^{(2)}(H_{0,1})\,,\ G^{(1)}(H_{0,1})\,\}
$$
for some meromorphic function $\mu$, or rather, meromorphic differential $\mu\,\dd z^2$.
Here we write $G^{(1)}(H_{0,1}) = \pi_1^{\p}(H_{0,1}^{(1)})$, a section of the Gauss bundle $G^{(1)}(\aal_1)$.
Similarly, we write $G^{(2)}(H_{0,1})$ for the component of the second derivative $H_{0,1}^{(2)}$ orthogonal to $(\aal_1)_{(1)} = \aal_1 \oplus
G^{(1)}(\aal_1)$; this is a section of the
second Gauss bundle $G^{(2)}(\al_1)$ of $\aal_1$ defined by \eqref{Gauss-i}.
Note that $\phi_2$ lies in a Grassmannian if and only if $\mu \equiv 0$, in which case $\phi_2 = G^{(1)}(\al_1)$.

This description is equivalent to that in \cite{Wo-Un}, but more economical.

\medskip

(iv) Let $\phi:M^2 \to \U(4)$ be a non-constant harmonic
 map of finite uniton number.  Then, precisely one the following occurs:

 \smallskip

($a$) $\phi$ has uniton number one and is given by a holomorphic map
$\phi:M^2 \to G_{d_1}(\C^3)$ where $d_1 = 1$, $2$ or $3$;

\smallskip

($b$) $\phi$ has uniton number two and is given by \eqref{fact} with the two unitons $\aal_1$, $\aal_2$ described as follows:

\indent\indent ($b_{12}$) \
$\aal_1$, $\aal_2$ have rank one and two respectively, with $\aal_1$ full,
and are given by \eqref{rank12};

\indent\indent ($b_{13}$) \
$\aal_1$, $\aal_2$ have rank one and three respectively, with $\aal_1$ full,
and are given by
$$
\aal_1 = \spn\{H_{0,1}\} \quad \text{and} \quad
\aal_2 = \spn\{H_{0,1} + \pi_1^{\p}H_{1,1}\,,\, \pi_1^{\p}H_{1,2}
	\,,\, \pi_1^{\p}H_{0,1}^{(1)}\,\}
$$
for some meromorphic vectors $H_{0,1}, H_{1,1}, H_{1,2}$\,;

\indent\indent ($b_{23}$) \
$\aal_1$, $\aal_2$ have rank two and three respectively, with $\aal_1$ full.
Because of the low dimensions, $\aal_1$ must have Gauss bundle
$G^{(1)} (\aal)$
 of rank one.  Now it is easily seen that the Gauss bundle is the image of the holomorphic bundle map $A^{\phi}_z|_{\aal_1}$ ; on choosing $H_{0,1}$ in the kernel of this map  we obtain
\begin{align*}
\aal_1 &= \spn\{H_{0,1}, H_{0,1}^{(1)}\},
\\
\aal_2 &= \spn\{H_{0,1} + \pi_1^{\p}H_{1,1}\,,\,
	H_{0,1}^{(1)} + \pi_1^{\p}H_{1,2}\,,\, \pi_1^{\p}H_{0,1}^{(1)}\,\},
\end{align*}
for some meromorphic vectors $H_{0,1}, H_{1,1}, H_{1,2}$\,;

\smallskip

(c) $\phi$ has uniton number three and is given by \eqref{fact} where the three unitons $\aal_1, \aal_2, \aal_3$ have rank one, two and three, respectively, and $\aal_1$ is full.  Then they are given by
\eqref{alpha1}---\eqref{alpha32} for a single column
$\bigl(H_{0,1}, H_{1,1}, H_{2,1}\bigr)^T$
of meromorphic vectors.

In contrast to previous descriptions in \cite{PiZa,Wo-Un}, these explicit formulae do not involve any unsolved $\ov{\pa}$-problems.
\end{exm}

For harmonic maps into Grassmannians, see the end of Section
\ref{sec:loop}.

\section{Loop groups and the Grassmannian model}\label{sec:loop}

A. \emph{Harmonic maps into $\U(n)$}

\medskip

Following \cite{PrSe} and \cite{Seg},
let $\Om^{\alg}\U(n)$ denote the \emph{algebraic loop group}
 consisting of all maps $\ga:S^1 \to \U(n)$ with $\ga(1) = I$
which are of the form
$$
\ga(\la) = \sum_{i=s}^t A_i \la^i
$$
for some integers $s,t$ and some $A_i \in \gl(n,\CC)$.
Let $\HH$ denote the Hilbert space $L^2(S^1,\C^n)$,
and let $\HH_+$ denote the closed subspace of elements of the form $\sum_{k \geq 0}\la^k a_k$ where $a_k \in \C^n$.  The loop group $\Om^{\alg}\U(n)$ acts on $\HH$ and the map $\ga \mapsto \ga(H_+)$ identifies
$\Om^{\alg}\U(n)$ with the
\emph{algebraic Grassmannian} consisting of all subspaces $W$ of $\HH$ such that $\la W \subseteq W$ and
$\la^s\HH_+ \subseteq W \subseteq \la^t\HH_+$ for some $s,t$.

Now, a polynomial extended solution $\Phi$ may be interpreted
 as a smooth map $\Phi:M^2 \to \Om^{\alg}\U(n)$.  Then
setting $W = \Phi(\HH_+)$ defines
 a map $W$ from $M^2$ to the algebraic Grassmannian,
and, it is easy to see that, if $\Phi$ has degree $r$,
the subspace satisfies $\la^r \HH_+ \subseteq W \subseteq \HH_+$, and so can be identified with the coset $W + \la^r\HH_+$ in the quotient space
$\HH_+/\la^r\HH_+$.
We thus obtain a map $W:M^2 \to G_*(\HH_+/\la^r\HH_+)$
into the Grassmannian of the finite-dimensional vector space
$\HH_+/\la^r\HH_+$\,.   Note that this vector space can
be canonically identified with $\C^{rn}$ via the isomorphism
\begin{equation} \label{iso}
(L_0,L_1,\ldots, L_{r-1}) \mapsto
L_0 + \la L_1 + \cdots \la^{r-1}L_{r-1}
	+ \la^r \HH_+ \,.
\end{equation}

Now, given a harmonic map $\phi:M^2 \to \U(n)$ of uniton number at
most $r$, let $\Phi$ be its unique type one extended solution.  Then setting
$W = \Phi(\HH_+)$ defines a holomorphic map
$W:M^2 \to \HH_+/\la^r \HH_+ \cong \C^{rn}$, or equivalently, a holomorphic subbundle $\WW$ of the trivial bundle $M^2 \times \HH_+/\la^r\HH_+ \cong \CC^{rn}$ satisfying
\begin{equation} \label{W-harm}
\la \WW_{(1)} \subseteq \WW \,;
\end{equation}
this, together with the holomorphicity, being equivalent to harmonicity of $\phi$, see for example \cite{Gu-bk}.   Here, for any $i \geq 0$, $\WW_{(i)}$ denotes the subbundle spanned by (local) sections of $\WW$ and their first $i$ derivatives (with respect to any complex coordinate on $M^2$); note that
 the condition \eqref{W-harm} incorporates the algebraic condition $\la\WW \subseteq \WW$.  We call $\WW$ the \emph{Grassmannian model} of $\phi$ (or $\Phi$).

M.\ Guest \cite{Gu-update} notes that holomorphic subbundles $\WW$ of $M^2 \times \HH_+/\la^r\HH_+$ satisfying
\eqref{W-harm} are given by taking an arbitrary holomorphic map $X:M^2 \to G_*(\HH_+/\la^r \HH_+) \cong G_*(\C^{rn})$, equivalently, holomorphic subbundle $\XX$ of $\CC^{rn}$, and setting $\WW$ equal to the coset
\begin{equation} \label{W}
\WW = \XX + \la \XX_{(1)} + \la^2 \XX_{(2)} + \cdots + \la^{r-1}  \XX_{(r-1)} + \la^r \HH_+ \,.
\end{equation}

Now, let $(\aal_1, \ldots,\aal_i)$ be a sequence of subbundles of $\CC^n$ for some $i \geq 1$.

For $0 \leq s \leq i$ write $S^{i}_{s}$ for the sum of all $i$-fold products
of the form $\Pi_i \cdots \Pi_1$  where exactly $s$ of the $\Pi_{\ell}$ are $\pi_{\ell}^{\p}$ and the other $i-s$ are $\pi_{\ell}$.  If $i=0$, set $S^{i}_{s} = I$.
Writing $\Phi_i = (\pi_1 + \la\pi_1^{\p})\cdots(\pi_i + \la\pi_i^{\p})$, we
see that $S^{i}_{s}$ is the coefficient of $\la^{-s}$ in
 $\Phi_i^{\;-1} = (\pi_i+ \la^{-1} \pi_i^{\p}) \cdots
	(\pi_1+ \la^{-1} \pi_1^{\p})$.

The $S^i_s$ are related to the $C^i_s$
defined by \eqref{Cis} in a simple way as follows; here
$\binom{i}{s}$ denotes the binomial coefficient
$i!/s!(i-s)!$\,.

\begin{lem} \label{lem:L-sum}
For all integers $k$ and $i$ with $0 \leq k \leq i$,
$$
\sum_{s=k}^i \binom{s}{k}S^i_s = C^i_k\,.
$$
\end{lem}

\begin{proof}
In the following proof, set $\binom{i}{s}$ and $S^i_s$ equal to $0$ if $s > i$ or $s <0$.  The result is clearly true for $i=0$ and $i=1$.   Suppose that it is true for $i-1$.  Then using
$$
\binom{i}{s} = \binom{i-1}{s-1} + \binom{i-1}{s}
\quad \text{and} \quad
S^i_s = \pi_i^{\p}S^{i-1}_{s-1} + \pi_i S^{i-1}_s
\,,
$$
we calculate
\begin{align*}
\sum_{s=k}^i \binom{s}{k} S^i_s
&=
\pi_i^{\p} \sum_{s=k}^i \binom{s}{k}S^{i-1}_{s-1}
+\pi_i \sum_{s=k}^i \binom{s}{k}S^{i-1}_{s} \\
&=
\pi_i^{\p} \sum_{s=k}^i \binom{s-1}{k-1}S^{i-1}_{s-1}
+\pi_i^{\p} \sum_{s=k+1}^i \binom{s-1}{k}S^{i-1}_{s-1}
+ \pi_i \sum_{s=k}^{i-1} \binom{s}{k}S^{i-1}_{s}\\
&=
\pi_i^{\p}C^{i-1}_{k-1} + \pi_i^{\p}C^{i-1}_{k}
+ \pi_i C^{i-1}_{k}\\
&=
\pi_i^{\p}C^{i-1}_{k-1} + C^{i-1}_k = C^i_k
\end{align*}
as required.
\end{proof}

For any $i \geq 0$ and meromorphic vectors $(H_0,H_1,\ldots, H_i)$, set
\begin{equation} \label{H-to-L}
L_i = \sum_{\ell=0}^i\binom{i}{\ell}H_{\ell}\,.
\end{equation}
As above, let $(\aal_1, \ldots,\aal_i)$ be a sequence of subbundles of $\CC^n$ for some $i \geq 1$.
Define $K_{\ell}^{(k)}$ \ $(0 \leq k \leq \ell \leq i-1)$ by \eqref{K-abstr},
and suppose that each $K_{\ell}^{(k)}$ is a section of\/
$\aal_{\ell+1}$. Then

\begin{lem} \label{lem:sum=Hi}
For all integers $i,j,k$ with $0 \leq k \leq j \leq i$,
$$
\sum_{s=j}^iS^i_s L_{s-j}^{(k)} =	\sum_{s=j}^i C^i_s H_{s-j}^{(k)} =
	\left\{ \begin{array}{cll} K_i^{(k)}	& \text{if} &j=k\,, \\
								0		& \text{if} &j>k\,.		
	\end{array} \right.
$$
\end{lem}

\begin{proof}
We start with the proof of the second identity. It clearly holds for $i =0$.
Suppose that it holds for some $i \geq 0$; we show that it holds for $i+1$.  For $k=j$, there is nothing to prove, as it is the definition of $K_i^{(k)}$, so we may assume that $k<j\leq i+1$.  Then
\begin{align*}
\sum_{s=j}^{i+1} C_s^{i+1} H_{s-j}^{(k)}
	&=\pi_{i+1}^{\p}\sum_{s=j}^{i+1} C_{s-1}^i H_{s-j}^{(k)}
		+\sum_{s=j}^{i+1} C_s^i H_{s-j}^{(k)}\\
	&=\pi_{i+1}^{\p}\sum_{s=j-1}^i C_s^i H_{s-(j-1)}^{(k)}
		+\sum_{s=j}^i C_s^i H_{s-j}^{(k)}
\end{align*}
Now, if $j=k+1$, the first term is equal to $\pi_{i+1}^{\p}(K_i^{(k)})$, which is zero since $K_i^{(k)}$ lies in $\aal_{i+1}$; on the other hand, if $j>k+1$, it is zero by the induction hypothesis.  The second term always vanishes by the induction hypothesis, completing the induction step.

As for the proof of the first identity, again it clearly holds for $i=0$, so suppose that it holds for some $i \geq 0$.  Then, for any $k\leq j\leq i+1$ we have
\begin{align*}
\sum_{s=j}^{i+1}S_s^{i+1}L_{s-j}^{(k)}
	&=\pi_{i+1}^{\p}\sum_{s=j}^{i+1}S_{s-1}^{i}L_{s-j}^{(k)}
		+\pi_{i+1}\sum_{s=j}^{i+1}S_s^{i}L_{s-j}^{(k)}\\
	&=\pi_{i+1}^{\p}\sum_{s=j-1}^{i}S_s^{i}L_{s-(j-1)}^{(k)}
		+\pi_{i+1}\sum_{s=j}^{i}S_s^{i}L_{s-j}^{(k)}\\
\intertext{(using the induction hypothesis on both terms)}
	&=\pi_{i+1}^{\p}\sum_{s=j-1}^{i}C_s^{i}H_{s-(j-1)}^{(k)}
		+\pi_{i+1}\sum_{s=j}^{i}C_s^{i}H_{s-j}^{(k)}\\
	&=\pi_{i+1}^{\p}\sum_{s=j}^{i+1}C_{s-1}^{i}H_{s-j}^{(k)}
		+(I-\pi_{i+1}^{\p})\sum_{s=j}^{i+1}C_s^{i}H_{s-j}^{(k)}\\
	&=\sum_{s=j}^{i+1}(\pi_{i+1}^{\p} C_{s-1}^{i}+C_s^{i})H_{s-j}^{(k)}
		-\pi_{i+1}^{\p}\sum_{s=j}^{i+1}C_s^{i}H_{s-j}^{(k)}\\
	&=\sum_{s=j}^{i+1}C_s^{i+1}H_{s-j}^{(k)}
		-\pi_{i+1}^{\p}\sum_{s=j}^{i}C_s^{i}H_{s-j}^{(k)}\,.
\end{align*}
Since, by the first part of the proof, $\sum_{s=j}^{i}C_s^{i}H_{s-j}^{(k)}$ is either $K_i^{(k)} \in \aal_{i+1}$ or zero, the last term is zero, and the induction step is complete.
\end{proof}

We thus have a second set of formulae for our unitons which are more adapted to the Grassmannian model, but not as easy to interpret geometrically.

In the following theorem, it is convenient to identify $\HH_+/\la^r \HH_+$ with $\CC^{rn}$ via the isomorphism \eqref{iso}.

\begin{thm} \label{thm:Grass-model}
Let $r \geq 1$, and let $\BB$ and $\XX$ be holomorphic subbundles of\/
$\CC^{rn}$ related by the linear isomorphism
$$
\BB \ni H = (H_0,H_1\ldots, H_{r-1}) \mapsto L = (L_0,L_1,\ldots,L_{r-1}) \in \XX
$$
given by \eqref{H-to-L}.
Let $\phi:M^2 \to \U(n)$ be the harmonic map given by $\BB$ according to Theorem \ref{thm:constr-B}, and let $W:M^2 \to G_*(\C^{rn})$ be the holomorphic map given by \eqref{W}.  Then $W$ is the Grassmannian model of $\phi$, i.e., the type one extended solution $\Phi$ of $\phi$ satisfies
$\Phi(\HH_+) = \WW$.
\end{thm}

This result follows from a purely algebraic result which constructs the factorization from $W$.  As above, we identify $\HH_+/\la^r \HH_+$ with $\C^{rn}$; let $P_i:\C^{rn} \to \C^n$ denote the canonical projection
$(L_0,L_1,\ldots,L_{r-1}) \mapsto L_i$.

\begin{prop} \label{prop:fact}
Let $r \geq 1$, and let $\WW$ be a subspace of $\HH_+/\la^r \HH_+$ with $\la \WW \subset \WW$.
Define subspaces $\aal_i$ \ $(i=1,\ldots, r)$ of\/ $\C^n$ by
\begin{equation} \label{W-to-al}
\aal_i = (\sum_{s=0}^{i-1} S^{i-1}_s P_s)\WW
\end{equation}
and define $\Phi \in \Omega^{\alg}(\U(n))$ by
\begin{equation} \label{Phi_r}
\Phi = (\pi_1+\la\pi_1^{\p}) \cdots	
			(\pi_r+\la\pi_r^{\p}) \,.
\end{equation}
Then
$\Phi(\HH_+) = \WW$.
\end{prop}

\begin{proof}
It is equivalent to show that $\Phi^{-1}(\WW) = \HH_+$\,.
For any $i = 1,\ldots, r$ write
$\Phi_i = (\pi_1+\la\pi_1^{\p}) \cdots (\pi_i+\la\pi_i^{\p})$, so that
$\Phi = \Phi_r$\,.  Then
$$
\Phi_i^{\;-1} = (\pi_i+\la^{-1}\pi_i^{\p})\cdots (\pi_1+\la^{-1}\pi_1^{\p})\,.
$$
We first show that $\Phi_i^{\;-1}(\WW) \subseteq \HH_+$; putting $r=i$ gives $\Phi^{-1}(\WW) \subseteq \HH_+$\,.  Let $L = L_0 + \la L_1 + \cdots + \la^{r-1}L_{r-1} + \la^r \HH_+$ be in $\WW$.  Then, for any $i <r$, the constant term of $\Phi_i^{\;-1}(L)$ is $K_i = \sum_{s=0}^i S^i_sL_s$ which lies in $\aal_{i+1}$.

We show that the expression for $\Phi_i^{\;-1}(L)$ contains no negative powers of $\la$.  This is clearly true for $i=0,1$; suppose that it is true for some value of $i \in \{0,1,\ldots, r-1\}$.  We have
$\Phi_{i+1}^{\;-1}(L)=(\pi_{i+1} \Phi_i^{\;-1} + \la^{-1}\pi_{i+1}^{\p}\Phi_i^{\;-1})L$.  By the
induction hypothesis, $\Phi_i^{\;-1}(L)$ lies in $\HH_+$,
so the only possible term with a negative power comes from
 $\la^{-1}\pi_{i+1}^{\p}$ applied to the constant term of $\Phi_i^{\;-1}(L)$.  But this gives $\la^{-1}\pi_{i+1}^{\p}(K_i)$ which is zero since $K_i$ lies in $\aal_{i+1}$.

We now show that $\Phi^{-1}(\WW) \supseteq \HH_+$.
Note that the constant term of $\Phi^{-1}(L) = \Phi_r^{\;-1}(L)$ is
$\sum_{s=0}^{r-1} S^r_s L_s + \pi_r^{\p}\cdots \pi_1^{\p}\HH_+$.  The first term is
$\pi_r K_{r-1} +\pi_r^{\p} \sum_{s=0}^{r-1} S^{r-1}_{s-1}L_s$ and, by the covering property, the second term equals $\pi_r^{\p} \HH_+$ so that
$\Phi^{-1}(L) = K_{r-1} + \aal_r^{\p}$.  It follows that the projection $P_0(\WW)$ contains $\aal_r + \aal_r^{\p} = \C^n$.
Since $\WW$ is closed under multiplication by $\la$\,,
so is $\Phi^{-1}(\WW)$.  It thus contains the whole of $\HH_+$\,.
\end{proof}

Note that, since any member $\ga$ of the algebraic complexified loop group gives such a subspace $W = \ga(\HH_+)$, then \eqref{W-to-al} and \eqref{Phi_r} give an \emph{explicit geometrical Iwasawa factorization} of that group.

\medskip

\emph{Proof of Theorem \ref{thm:Grass-model}.}

The formula in Lemma \ref{lem:sum=Hi} for $j=k=0$ tell us
$$
\aal_i^{(0)}
= \bigl\{\sum_{s=0}^{i-1} \bigl( S_s^{i-1}P_s \bigr)L : L \in X \bigr\}\,.
$$
Since $\la L^{(1)} = (0,L_0^{(1)},L_1^{(1)},\ldots, L_{r-2}^{(1)})$,
putting $j=k=1$ gives
$$
\aal_i^{(1)}
= \bigl\{\sum_{s=0}^{i-1} \bigl( S_s^{i-1}P_s \bigr)\la L^{(1)} : L \in X \bigr\}\,.
$$
and similarly for the other $\aal_i^{(k)}$.

Further, since $\la L = (0,L_0,L_1,\ldots, L_{r-2})$, putting $k=0$ and $j=1$ gives
$$
\bigl\{\sum_{s=0}^{i-1} \bigl( S_s^{i-1}P_s \bigr)\la L : L \in X \bigr\}
= \bigl\{\sum_{s=1}^{i-1} S_s^{i-1}L_{s-1}: L \in X \bigr\}
 = \vec{0}\,.
$$
Combining these and using \eqref{W}, we see that \eqref{W-to-al} holds, and so the calculations of Proposition \eqref{prop:fact} show that
$\Phi(\HH_+) = \WW$, as desired.
\qed

\begin{exm}
The Grassmannian models corresponding to the harmonic maps
in Example \ref{exm:low-dim} (iii) and (iv) are
given by \eqref{W} with $X$ as follows:

\smallskip
\emph{Maps into $\U(3)$.}

\indent\indent\emph{Case $(a)$}: $X = \spn\{L_{0,j}: 1\leq j \leq d_1\}$,

\indent\indent\emph{Case $(b)$}: $X = \spn\{L_{0,1} + \la L_{1,1}$\}.

\medskip

\emph{Maps into $\U(4)$.}

\indent\indent\emph{Case $(a)$}: $X = \spn\{L_{0,j}: 1\leq j \leq d_1\}$;

\indent\indent\emph{Case $(b_{12})$}: $X = \spn\{L_{0,1} + \la L_{1,1}\}$,

\indent\indent\emph{Case $(b_{13})$}: $X = \spn\{L_{0,1} + \la L_{1,1}\,,\, \la L_{1,2}\}$,

\indent\indent\emph{Case $(b_{23})$}: $X = \spn\{L_{0,1} + \la L_{1,1}\,,\,
	L_{0,1}^{(1)} + \la L_{1,2}\}$;

\indent\indent\emph{Case $(c)$}:  $X = \spn\{L_{0,1} + \la L_{1,1} + \la^2 L_{2,1}\}$.

\smallskip

These are exactly the cases found by Guest \cite{Gu-update} as the Grassmannian models corresponding to the harmonic maps of \cite{BuGu}.
\end{exm}

\bigskip

B. \emph{Harmonic maps into Grassmannians}

\medskip

Using the Grassmannian model, we can give construct harmonic maps into Grassmannians.  Let $\Phi_{\la}$ be the type one
 extended solution of some harmonic map from a surface to
$\U(n)$; we have the following characterization due to Uhlenbeck \cite{Uhl}.

\begin{lem} \label{lem:sym}
Let\/ $\phi:M^2 \to \U(n)$ be a harmonic map and let\/ $\Phi_{\la}$ be its type one extended solution, so that\/
$\phi =Q \Phi_{-1}$ for some $Q \in \U(n)$.  Then $\phi$ maps into a Grassmannian if and only if\/
$Q^2 = I$ and
\begin{equation} \label{twist-Q}
\Phi_{\la} = Q \Phi_{-\la}\Phi_{-1}^{\;-1}Q.
\end{equation}
\end{lem}

\begin{proof}
If $Q^2 = I$ and \eqref{twist-Q} holds, then, putting $\la=-1$ shows that $\phi$ is Grassmannian.  Conversely, given a harmonic map $\phi$ into a Grassmannian, Uhlenbeck \cite[page 25]{Uhl} first constructs an extended
solution $E_{\la}$ with $E_{-1} = \phi$. Then, in her Lemma 15.1, she shows that
the type one extended solution of $\phi$ satisfies $\Phi_{\la} = Q(\la)^{-1}E_{\la}$ for some $Q(\la)$
with $Q(\la) = Q(-\la)Q(-1)$, so that
$\Phi_{-1} = Q(-1)^{-1}E_{-1} = Q\phi$
where $Q = Q(-1) = Q(-1)^{-1}$; she further shows that \eqref{twist-Q} holds for this same $Q$.
\end{proof}

Note that (i) the condition $Q^2 = I$ actually follows from \eqref{twist-Q} by putting $\la=1$, it says that $Q$ lies in a Grassmannian;
(ii) since $\Phi$ is uniquely determined by $\phi$, so is $Q$. However, if we are given only $\Phi$, there may be more than one $Q$ satisfying \eqref{twist-Q}, and so more than one harmonic
map $\phi$ into a Grassmannian with the same type one extended solution $\Phi$\,; this problem will be studied elsewhere.

The lemma implies the following criterion in terms of the Grassmannian model.

\begin{prop}
 \label{prop:Grass-Q}
Let $\Phi$ be the type one extended solution of a harmonic map of finite uniton number, say\/ $r$, and let
$W = \Phi\HH_+ \in \HH_+/\la^r\HH_+$ be the corresponding Grassmannian solution.
Then $\Phi$ is the extended solution of a harmonic map into a Grassmannian if and only if there exists
$Q \in \U(n)$ with $Q^2 = I$ such that
$W = \Phi \HH_+$ is closed under the involution
$\nu_Q:\HH_+/\la^r\HH_+ \to \HH_+/\la^r\HH_+$ defined by
$$
L(\la) = \sum_{i=0}^{r-1} L_i \la^i
	\mapsto Q L(-\la) = \sum_{i=0}^{r-1} (-1)^i Q L_i \la^i,
$$
i.e.,
\begin{equation} \label{W-Q}
W_{\la} = Q W_{-\la} \quad (\la \in \C) \,.
\end{equation}
\end{prop}

\begin{proof}
We shall prove that \eqref{twist-Q} holds for some $Q$ if and only if \eqref{W-Q} holds.  The result then follows from Lemma \ref{lem:sym}.

Since $\Phi_{-1}^{\;-1} Q$ acts isomorphically on $\HH_+$, the subspace $W$ satisfies
$$
W_{\la} = \Phi_{\la}(\HH_+) = Q\Phi_{-\la}\Phi_{-1}^{\;-1} Q \HH_+
 = Q\Phi_{-\la}\HH_+ = Q W_{-\la}
$$
 for all $\la \in S^1$, and so, by real analyticity, for all $\la \in \C$.

Conversely, suppose that $W$ is closed under the involution
$\nu_Q$.  Then
$\Phi_{\la}\HH_+ = Q\Phi_{-\la}\HH_+$.
Hence $(Q\Phi_{-\la})^{-1}\Phi_{\la}$ stabilizes $\HH_+$, and is so a constant in $\U(n)$.  Putting $\la= \pm 1$ shows that this constant is $\Phi_{-1}^{\;-1}Q = Q\Phi_{-1}$, and we obtain \eqref{twist-Q}.
\end{proof}

To express this explicitly note that any $Q \in \U(n)$ with $Q^2 = I$ is of the form $Q = \pi_A - \pi_A^{\p}$ for some unique subspace $A$ (in fact, $A$ and $A^{\p}$ are the $\pm 1$-eigenspaces of $Q$). Then say that a polynomial
$L \in \HH_+/\la^r \HH_+$ is \emph{$Q$-adapted} if its
coefficients have image alternately in $A$ and $A^{\p}$, i.e.,
$L(\la) = \sum_{i=0}^{r-1} L_i \la^i$ and \emph{either}

{\rm ($+$)} $L_i$ has image in $A$ for $i$ even, and in $A^{\p}$ for $i$ odd, \emph{or}

{\rm ($-$)} $L_i$ has image in $A^{\p}$ for $i$ even, and in $A$ for $i$ odd;

\noindent
equivalently, $L$ lies in \emph{either} the $(+1)$- \emph{or} the
$(-1)$-eigenspace of $\nu_Q$.

\begin{cor}
$\Phi$ is the extended solution of a harmonic map into a Grassmannian if and only if\/ $W$ has a spanning set consisting of\/ $Q$-adapted polynomials, or,
equivalently, $W$ is given by \eqref{W} for some $X$ which has a
spanning set consisting of\/ $Q$-adapted polynomials.
\end{cor}

\begin{proof}
If $W$ has a basis with each element of the type $(+)$ or $(-)$ it is clear that $W$ is closed under $\nu_Q$.  Conversely,
if $W$ is closed under $\nu_Q$, then given an arbitrary spanning set $\{e_i\}$ for $W$, the set $\{e_i \pm \nu_Q e_i\}$ is also a spanning set, but now with each element of the type $(+)$
or $(-)$.  A subset of these will give a basis.

For the last statement, note that if $L$ is of type $(+)$ (resp.\ $(-)$), then $\la L$ and $\la L^{(1)}$ are of type $(-)$ (resp.\ $(+)$); hence if $X$ has a $Q$-adapted spanning set, so does $W$.
\end{proof}

Note that \emph{when $Q = I$, a polynomial is $Q$-adapted if
and only if it is even or odd, i.e, has coefficients of all odd or all even powers of $\la$ equal to zero}.

\medskip

Using \eqref{H-to-L} as in Theorem \ref{thm:Grass-model}, we can translate these condition on $X$ into conditions on $B$, i.e, on the original data $(H_{i,j})$.  A geometrical treatment of this will be given elsewhere; we content ourselves here with some simple examples where the geometry is transparent.

\begin{exm} \label{exm:Grass-non-isotr}

{\rm (i)}
With $r=3$, let $L_{0,1}$ and $L_{2,1}$ be arbitrary meromorphic vectors,
and let $X$ be spanned by the even polynomial
$L_{0,1} + \la^2 L_{2,1}$.  Then $W$ is spanned by this and by the odd polynomial $\la L_{0,1}^{(1)}$; both polynomials are
$Q$-invariant with $Q=I$.

Via \eqref{H-to-L}, $X$ corresponds to a single column $(H_{0,1},0,H_{2,1})^T$ of data where $H_{0,1} = L_{0,1}$ and
$H_{2,1} =L_{2,1}-L_{0,1}$, but it is easy to see that taking $H_{0,1} = L_{0,1}$ and $H_{2,1} = L_{2,1}$ gives the same unitons and so the same harmonic map.

For $n=4$, this is the special case of Example \ref{exm:low-dim}(iv)(c) given by $H_{1,1} = 0$.   We see that
$\aal_1 = \spn\{H_{0,1}\}$,
$\aal_2 = (\aal_1)_{(1)} = \spn\{H_{0,1}, H_{0,1}^{(1)}\}$ and
$\aal_3 = \spn\{H_{0,1} + \pi_2^{\p}H_{2,1}, \pi_1^{\p}H_{0,1}^{(1)}, \pi_2^{\p}H^{(2)}_{0,1}\}$, so that
our formula \eqref{fact} gives the harmonic map $\phi: M^2 \to G_2(\C^4)$ given by
$$
\pphi =\spn\{H_{0,1} + \pi_2^{\p}H_{2,1}, \pi_2^{\p}H_{0,1}^{(2)}\}\,.
$$
Assuming that
$\aal_1$ is full, $\phi$ has uniton number three.

\medskip

{\rm (ii)}
With $r=2$, let $Q = \pi_A - \pi_A^{\p}$ where $A$ is a proper subspace of $\C^n$. Let $L_{i,j}$ $(i=0,1,\,j=1,2)$ be meromorphic vectors where $L_{0,1}$ and $L_{1,2}$ have values in $A$ but $L_{0,2}$ and $L_{1,1}$ have values in $A^{\p}$.
Let $X$ be spanned by
$L_{0,j} + \la L_{1,j}$ \ $(j=1,2)$.
Then $W$ is spanned by these and by $\la L_{0,j}^{(1)}$ \ $(j=1,2)$, and all four polynomials are $Q$-adapted.
Note that \eqref{H-to-L} gives $H_{0,j} = L_{0,j}$ and $H_{1,j} = L_{1,j}-L_{0,j}$, but we obtain the same unitons, and so the same harmonic map, from the data
$H_{i,j} = L_{i,j}$\,.

Then the formula \eqref{W-to-al} gives the first uniton
$\aal_1$ spanned by $L_{0,1}=H_{0,1}$ and $L_{0,2}=H_{0,2}$; we can choose this data such that $\aal_1$ is full --- it suffices to have no linear relation between their components.

{}From \eqref{W-to-al}, the
second uniton $\aal_2$ is the rank 4 subbundle of $\CC^n$ spanned by $K_{1,j}= \pi_1 L_{0,j} + \pi_1^{\p} L_{1,j} = H_{0,j} + \pi_1^{\p} H_{1,j}$ and
$K_{1,j}^{(1)} = \pi_1^{\p} L_{0,j}^{(1)} = \pi_1^{\p} H_{0,j}^{(1)} $ \ $(j=1,2)$,
and $\phi = Q(\pi_1 - \pi_1^{\p})(\pi_2 - \pi_2^{\p})$ is a harmonic map of uniton number two into a Grassmannian.
In fact, $\phi_1 = Q(\pi_1-\pi_1^{\p})$ is the harmonic map into $G_*(\C^n)$ given by
$\pphi_1 = \spn\{H_{0,1}\} \oplus \bigl(\spn\{H_{0,2}\}^{\p}\cap A^{\p}\bigr)$, so that
$$
\pphi = \spn\{K_{1,1}, K_{1,2}^{(1)}\}\oplus \spn\{K_{1,2}, K_{1,1}^{(1)}\}^{\p}\! \cap \pphi_1^{\p}.
$$
However, it is easily checked that $\Phi_{-1} = (\pi_1 - \pi_1^{\p})(\pi_2 - \pi_2^{\p})$ does \emph{not} have image in a Grassmannian.
\end{exm}

\medskip

We finish with the important class of \emph{$S^1$-invariant} harmonic maps.

\begin{exm} \label{exm:isotropic}
Suppose that $X$ is spanned by monomials of the form
$\la^k L_k$; equivalently,
the data $(H_{i,j})$ is in the `diagonal' form
\begin{equation*}
\left[
\begin{array}{ccccccccccc}
H_{0,1} & \cdots & H_{0,d_1} & 0 &\cdots & 0
	&0 & \cdots & 0 & 0 & \cdots \\
	0  	& \cdots & 0 & H_{1,d_1+1} &
	\cdots & H_{1,d_2} & 0 & \cdots & 0 & 0 & \cdots \\
	0	& \cdots & 0 & 0 & \cdots & 0 &
	H_{2,d_2+1} & \cdots & H_{2,d_3} & 0 & \cdots \\
\vdots & \vdots & \vdots & \vdots & \vdots & \vdots
	& \vdots & \vdots & \vdots & \vdots & \ddots \\	
\end{array} \right]
\end{equation*}
then the unitons $\aal_i$ are \emph{nested}, i.e.,
$\aal_i \subseteq \aal_{i+1}$.
The projections in \eqref{alpha-H} are unnecessary and
$$
\aal_{i+1}
= \spn\{H_{i-\ell,k}^{(\ell)}:1\leq k \leq d_{i+1}^{(0)},0\leq \ell \leq i\}\,.
$$
The harmonic map
$\phi = (\pi_1-\pi_1^{\p}) \cdots (\pi_r-\pi_r^{\p})$ has image in a
Grassmannian and is determined by the following formula which gives
 $\pphi$ if $r$ is odd or $\pphi^{\p}$ if $r$ is even.
\begin{equation} \label{phi-nested}
\pphi(^{\p}) =
\sum_{k=0}^{[(r-1)/2]}\aal_{r-1-2k}^{\p} \cap \aal_{r-2k} =
\aal_{r-1}^{\p} \cap \aal_r \oplus \aal_{r-3}^{\p} \cap \aal_{r-2} \oplus \cdots\,,
\end{equation}
where we set $\aal_0$ equal to the zero subbundle.  We thus obtain a harmonic map into a Grassmannian which is invariant under a natural $S^1$ action, see \cite{DoEs,Uhl}.
\end{exm}


\begin{thebibliography}{}
\bibitem{BuGu} F.~E.\ Burstall and M.~A.\ Guest,
\textit{Harmonic two-spheres in compact symmetric spaces, revisited},
Math.\ Ann.\ \textbf{309} (1997) 541--572.

\bibitem{BuWo} F.~E.\ Burstall and J.~C.\ Wood,
\textit{The  construction  of  harmonic  maps
into  complex Grassmannians},
J. Diff. Geom. \textbf{23} (1986) 255--298.

\bibitem{DaTe} B.\ Dai and C.-L.\ Terng,
\textit{B\"acklund transformations, Ward solitons, and unitons},
J. Differential Geom. \textbf{75} (2007) 57--108.

\bibitem{DoEs} J.\ Dorfmeister and J.-H.\ Eschenburg,
\textit{Pluriharmonic maps, loop groups and twistor theory},
Ann. Global Anal. Geom. \textbf{24} (2003) 301--321.

\bibitem{Gu-bk} M.~A.\ Guest,
\textit{Harmonic maps, loop groups, and integrable systems}, London Mathematical Society Student Texts, 38,
Cambridge University Press, Cambridge, 1997.

\bibitem{Gu-update} M.~A.\ Guest,
\textit{An update on harmonic maps of finite uniton number, via the zero curvature equation},
Integrable systems, topology, and physics (Tokyo, 2000),  85--113,
Contemp.\ Math. \textbf{309}, Amer. Math. Soc., Providence, RI, 2002.

\bibitem{HeSh} Q.\ He and Y.~B.\ Shen,
\textit{Explicit construction for harmonic surfaces in ${\rm U}(N)$ via adding unitons},
Chinese Ann. Math. Ser. B  \textbf{25}  (2004) 119--128.

\bibitem{KoMa} J.~L.\ Koszul and B.\ Malgrange,
\textit{Sur certaines structures fibr\'ees complexes},
Arch. Math. \textbf{9} (1958) 102--109.

\bibitem{PiZa} B.\ Piette and W.~J.\ Zakrzewski,
\textit{General solutions of the ${\rm U}(3)$ and ${\rm U}(4)$ chiral $\sigma$ models in two dimensions},
Nuclear Phys. B \textbf{300} (1988) 207--222.

\bibitem{PrSe} A.\ Pressley and G.\ Segal,
\textit{Loop groups},
Oxford Mathematical Monographs, Oxford Science Publications, The Clarendon Press, Oxford University Press, Oxford, 1986.

\bibitem{Seg} G.\ Segal,
\textit{Loop groups and harmonic maps},
Advances in homotopy theory (Cortona, 1988), 153--164,
London Math. Soc. Lecture Note Ser., 139,
Cambridge Univ. Press, Cambridge, 1989.

\bibitem{SvWo} M.\ Svensson and J.~C.\ Wood
\textit{Filtrations, factorizations and explicit formulae for harmonic maps},
preprint (2009).

\bibitem{Uhl} K.\ Uhlenbeck,
\textit{Harmonic maps into Lie groups: classical solutions of the chiral model},
J. Differential Geom. \textbf{30} (1989) 1--50.

\bibitem{Wol} J.~G.\ Wolfson,
\textit{Harmonic sequences and harmonic maps of surfaces into complex Grassmann manifolds},
J. Differential Geom. \textbf{27} (1988) 161-178.

\bibitem{Wo-Un} J.~C.\ Wood,
\textit{Explicit construction and parametrization of harmonic
two-spheres in the unitary group},
Proc. London Math. Soc. (3) \textbf{58} (1989) 608--624.

\end{thebibliography}
\end{document}